\documentclass[11pt]{article}
\usepackage[english]{babel}
\usepackage[T1]{fontenc}
\usepackage[utf8]{inputenc}
\usepackage{amsbsy}
\usepackage{amsmath}
\usepackage{amssymb}
\usepackage{amsfonts}
\usepackage{amsthm}
\usepackage{bm}
\usepackage{graphicx}
\usepackage[active]{srcltx}
\usepackage{upgreek}
\usepackage{xcolor}
\usepackage[all]{xy}

\usepackage{subcaption}

\usepackage{tkz-tab}

\newcommand{\C}{\mathbb{C}}

\newcommand{\N}{\mathbb{N}}

\newcommand{\R}{\mathbb{R}}

\newcommand{\Z}{\mathbb{Z}}

\newcommand{\gp}{\mathfrak{p}}
\newcommand{\gq}{\mathfrak{q}}

\newcommand{\gs}{\mathfrak{s}}

\newcommand{\gv}{\mathfrak{v}}

\newcommand{\grad}{\nabla}

\newcommand{\dive}{\mathrm{div}}

\newcommand{\energyset}{\mathcal{X}^1(\R)}
\newcommand{\Nenergyset}{\mathcal{NX}^1(\R)}

\newcommand{\petitodeunnii}{\underset{n\rightarrow +\ii}{o(1)}}

\newcommand{\ntend}{\underset{n\rightarrow +\ii}{\longrightarrow}}

\newcommand{\ii}{\infty}

\newcommand{\Emin}{E_{\mathrm{min}}}

\newcommand{\ntendf}{\underset{n\rightarrow +\ii}{\rightharpoonup}}

\newcommand{\sech}{\mathrm{sech}}

\newcommand{\loc}{\mathrm{loc}}

\newcommand{\normLii}[1]{\Vert #1\Vert_{\ii}}
\newcommand{\normLdeux}[1]{\Vert #1\Vert_{L^2}}
\newcommand{\normLun}[1]{\Vert #1\Vert_{L^1}}

\newtheorem{thm}{Theorem}[section]

\newtheorem{cor}[thm]{Corollary}
\newtheorem{lem}[thm]{Lemma}
\newtheorem{prop}[thm]{Proposition}

\newtheorem*{thm*}{Theorem}

\newtheorem{rem}[thm]{Remark}
\newtheorem{de}[thm]{Definition}

\theoremstyle{remark}

\title{Minimizing travelling waves for the one-dimensional nonlinear Schrödinger equations with non-zero condition at infinity}

\begin{document}
\author{
\renewcommand{\thefootnote}{\arabic{footnote}}
Jordan Berthoumieu\footnotemark[1]}
\footnotetext[1]{CY Cergy Paris Universit\'e, Laboratoire Analyse, G\'eom\'etrie, Mod\'elisation, F-95302 Cergy-Pontoise, France. E-mail: {\tt jordan.berthoumieu@cyu.fr}}
\maketitle

\begin{abstract}
This paper deals with the existence of travelling wave solutions for a general one-dimensional nonlinear Schrödinger equation. We construct these solutions by minimizing the energy under the constraint of fixed momentum. We also prove that the family of minimizers is stable. Our method is based on recent articles about the orbital stability for the classical and nonlocal Gross-Pitaevskii equations~\cite{BetGrSa2, deLaMen1}. It relies on a concentration-compactness theorem, which provides some compactness for the minimizing sequences and thus the convergence (up to a subsequence) towards a travelling wave solution.
\end{abstract}

\section{Introduction}\label{Introduction}
We are interested in the defocusing nonlinear Schrödinger equation

\begin{equation}\label{NLS}\tag{NLS}
    i\partial_t \Psi +\Delta \Psi +\Psi f(|\Psi|^2)=0\quad\text{on }\R\times\R.
\end{equation}

This equation appears as a relevant model in condensed matter physics. In particular, it is relevant in the context of the Bose-Einstein condensation or superfluidity (see~\cite{AbiHuMeNoPhTu,GinzPit1,Pitaevs1,Gross1,Coste1}) and in nonlinear optics (see~\cite{KivsLut1}), when the natural condition at infinity is 

\begin{equation}\label{nonvanishing condition à l'infini}
    |\Psi(t,x)|\underset{|x|\rightarrow +\ii}{\longrightarrow}1.
\end{equation}

This condition differs from the case of null condition at infinity, in the sense that the dispersion relation is different. In \eqref{NLS}, the function $f$ can be taken equal to $f(\rho)=1-\rho$. $\!$We obtain the Gross-Pitaevskii equation, but we can also take many other functions that provide possible alternative behaviours as enumerated by D. Chiron in~\cite{Chiron7}. In order to stay close to the behaviour of the Gross-Pitaevskii equation and to remain consistent with the nonvanishing condition~\eqref{nonvanishing condition à l'infini}, we shall assume that $f$ satisfies $f(1)=0$.

The equation is Hamiltonian. Its hamiltonian, the generalized Ginzburg-Landau energy, is~given~by

\begin{equation}\label{définition de l'energie de Ginzburg Landau}
    E(\Psi):=\int_\R e(\Psi):=E_k(\Psi)+E_p(\Psi):=\dfrac{1}{2}\int_\R |\partial_x \Psi|^2+ \dfrac{1}{2}\int_\R F(|\Psi|^2),
\end{equation}

with

\begin{equation}\label{definition de F}
    F(\rho):=\int_\rho^1 f(r)dr.
\end{equation}

We also introduce the (renormalized) momentum defined for a non-vanishing function $\Psi$, by the formula

\begin{equation*}
    p(\Psi)=\dfrac{1}{2}\int_\R \langle \Psi,i\partial_x\Psi \rangle_\C\Big(1-\dfrac{1}{|\Psi|^2}\Big),
\end{equation*}

where $\langle , \rangle_\C$ denotes the usual real scalar product defined by $\langle x,y\rangle_\C = \mathrm{Re}(x\overline{y}), \forall x,y\in\C$.

Those quantities are defined and conserved at least formally. In the sequel, we shall restrict the study to the case where $F$ is a nonnegative function, and we will focus on the Hamiltonian framework in which all the functions have finite energy.

If $\Psi$ does not vanish, we can apply the Madelung transform $\Psi=\sqrt{\rho}e^{i\varphi}$ where $\rho$ and $\varphi$ are as smooth as $f$ is. These variables satisfy the hydrodynamical form of the equation 

\begin{equation}\label{forme hydrodynamique de NLS}
    \left\{
\begin{array}{l}
    \partial_t \rho +\dive(\rho v)=0, \\
    \partial_t v +v\grad v = \grad\Big( \dfrac{|\grad \rho|^2}{2\rho^2}-\dfrac{\Delta \rho}{\rho}\Big)+2 f'(\rho)\grad \rho,
\end{array}
\right.
\end{equation}

where $v=2\grad \varphi$. By linearizing this system around the trivial solution $(\rho,v)=(1,0)$, this linearized system reduces, in the long wave approximation, to the free wave equation, with the sound speed

\begin{equation}\label{relation entre vitesse c_s et f'(1)}
c_s = \sqrt{-2 f'(1)},
\end{equation}

when the additional condition

\begin{equation}\label{f'(1) <0}
    f'(1)<0,
\end{equation}
is fulfilled, which we will assume throughout the paper.

We focus on the one-dimensional travelling waves. They are solutions of \eqref{NLS} of the form \begin{equation*}
    \Psi(t,x)=u(x+ct)\quad\text{for }(t,x)\in\R\times\R,
\end{equation*}

where $c\in\R$ is the speed of the travelling wave\footnote{We can restrict to nonnegative speeds, noticing that if $\Psi$ is a  travelling wave solution of \eqref{NLS}, then $\overline{\Psi}$ is a travelling wave solution of speed $-c$.}. Their profile $u$ is solution of the equation 

\begin{equation}\label{TWC}\tag{$TW_{c}$}
    ic u' +  u'' + u f(|u|^2)=0.
\end{equation}

In the sequel, we will label by $\mathfrak{v}_c$ non constant travelling waves with speed $c$. In particular, $\mathfrak{v}_0$ will be a stationary solution of the Schrödinger equation.$\,\!$ In the case of the Gross-Pitaevskii equation, non constant finite energy travelling waves exist for any speed $c\in (-\sqrt{2},\sqrt{2})$ and they are unique, up to a translation and a constant phase shift. Their shape was explicitly computed in the physical literature (see~\cite{BetGrSa2} for a rigorous description). For a general nonlinearity $f$, Z. Lin gave in~\cite{LinZhiw1} a sufficient and necessary condition for their existence and uniqueness. This condition is related to a general result concerning ordinary differential equations due to H. Berestycki and P.-L. Lions in Theorem~5 in~\cite{BereLions1}.

\begin{thm}[\cite{Chiron7,LinZhiw1}]\label{existence et unicité des travelling wave si condition du zéro}
Let $c\in[0,c_s]$. We assume that there exists
$\xi_c\in\! (-\ii,1]\setminus \{0\}$ such that $\mathcal{N}_c(\xi_c)=0$, $\mathcal{N}_c(\xi) <0$ on $\left\{
\begin{array}{l}
    ]0,\xi_c[,\text{ and }\mathcal{N}'_c(\xi_c)>0\text{ if }0 < \xi_c\leq 1, \\
    ]\xi_c,0[,\text{ and }\mathcal{N}'_c(\xi_c)< 0\text{ if } \xi_c < 0,
\end{array}
\right.$ where $\mathcal{N}_c(\xi)=c^2 \xi^2 -4(1-\xi)F(1-\xi)$.
Then there exists a unique non constant solution $\mathfrak{v}_c$ of \eqref{TWC}, up to a translation and a constant phase shift, that satisfies $|\gv_c(x)|\underset{|x|\rightarrow +\ii}{\longrightarrow}1$. The other solutions are the constant functions of modulus one.
\end{thm}

\begin{rem}
This proof is based on applying the arguments in~\cite{BereLions1} to the equation satisfied by $\eta:=1-|u|^2$. This provides a real-valued radial and decreasing solution $\eta_c$ from which the existence of the complex valued solution $\gv_c$ is deduced.
\end{rem}

\begin{rem}
    Theorem~\ref{existence et unicité des travelling wave si condition du zéro} also provides exponential decay to 0 at infinity of $1-|\gv_c|^2,\gv_c'$ and $\gv_c''$ (when assumption~\eqref{hypothèse de croissance sur F majorant} below is satisfied), so that such travelling waves have finite energy. 
\end{rem}

\begin{rem}\label{remarque non existence de tw supersonique sonique Maris}
The study necessarily reduces to the case $c\in [0,c_s)$, provided that $f$ is smooth enough. Indeed, we anticipate the fact that a sufficient condition for the orbital stability is that $f''(1)+3f'(1)$ is not zero. Going to Theorem 5.1 in~\cite{Maris6} and to the remark just after, we can claim in this case that there is no sonic or supersonic non constant travelling wave of finite energy.
\end{rem}

The energy $E$ is formally conserved by the flow of \eqref{NLS}, but the proof of this conservation first requires to give a proper sense to this energy. This is why we introduce the energy sets

\begin{equation*}
    \energyset:=\big\lbrace u\in H^1_{\loc}(\R) \big| u' \in L^2 (\R), F(|u|^2) \in L^1(\R) \big\rbrace 
\end{equation*}

and 

\begin{equation*}
    \Nenergyset:=\big\lbrace u\in \energyset \big| \inf_\R |u| >0 \big\rbrace .
\end{equation*}

Note that this energy set is exactly the same one as in the Gross-Pitaevskii case under the assumptions of Theorem~\ref{theoreme existence de travelling wave} below (see Remark~\ref{remarque: conséquence de H1 et H2, estimée des normes}). This is useful for addressing the Cauchy problem for~\eqref{NLS}. Indeed, a preliminary step for dealing with orbital stability is the well-posedness of the Cauchy problem $\eqref{NLS}$ with the nonvanishing condition at infinity \eqref{nonvanishing condition à l'infini}. The Cauchy problem for the Gross-Pitaevskii equation was solved in the energy space $\mathcal{X}^1(\R^N)$ by P. Gérard in~\cite{Gerard1} for $N\in\{2,3\}$ and for $N=4$ by R. Killip, T. Oh, O. Pocovnicu and M. Visan in~\cite{KiOhPoVi}. P. Zhidkov showed in~\cite{Zhidkov0} the local well-posedness of the Cauchy problem \eqref{NLS} on the space $\mathcal{Z}^k(\R):=\lbrace u \in L^{\ii}(\R)|\grad u \in H^{k-1}(\R)\rbrace$ (for $k\geq 1$). See also the article of C. Gallo~\cite{Gallo3} for the same result on $\mathcal{Z}^k(\R^N)$, provided that $k >\frac{N}{2}$ and for the rigorous justification of the energy and momentum conservation (if $k=N=1$ or $2$)\footnote{Only if $k=N=1$ for the momentum.}. The energy conservation combined to the fact that the equation is defocusing yields the global well-posedness. Under suitable conditions on $f$, C. Gallo finally showed in~\cite{Gallo1} that for $N\leq 4$, the Cauchy problem~\eqref{NLS} is globally well-posed in $u_0 + H^1(\R^N)$, provided that the initial condition $u_0$ is in the energy space defined below. Observing that $\energyset + H^1(\R) \subset \energyset$, we get a proper framework for well-posedness before addressing the question of orbital stability.

\begin{thm}[Theorem 1.2 in C. Gallo~\cite{Gallo1}]\label{global well-posedness of cauchy problem}
Let $u_0\in \energyset$. Take $f$ in $\mathcal{C}^2(\R)$ satisfying \eqref{1ere hypothèse de croissance sur F minorant intermediaire} below. In addition, assume that there exist $\alpha_1\geq 1$ and $C_0 >0$ such that for all $\rho\geq 1$,

\begin{equation}\label{theoreme de gallo condition de croissance sur f''}
|f''(\rho)|\leq \dfrac{C_0}{\rho^{3-\alpha_1}}.
\end{equation}

If $\alpha_1 > \frac{3}{2}$, assume moreover that there exists $\alpha_2\in [ \alpha_1-\frac{1}{2},\alpha_1]$ such that for $\rho\geq 2$, $C_0 \rho^{\alpha_2} \leq F(\rho)$.\\
There exists a unique function $w \in \mathcal{C}^0\big(\R,H^1(\R)\big)$ such that $u:=u_0+w$ solves \eqref{NLS}. Moreover, the solution depends continuously on the initial condition, and the energy $E$ and the momentum $p$ are conserved by the flow.
\end{thm}

We can also characterize the travelling waves by noticing that the equation~\eqref{TWC} can be formally written $\nabla E(\gv_c)=c\nabla p(\gv_c)$. This is the Euler-Lagrange equation associated with the minimization of the energy when the momentum is fixed, where $c$ appears as a Lagrange multiplier and $\gv_c$ a solution to this variational problem. In this context, we consider for $\gp\in\R$,

\begin{equation}\label{définition de la minimization curve Emin}
    E_{\min}(\mathfrak{p}):=\inf\big\lbrace E(v)\big| v\in\mathcal{N}\energyset, p(v)=\mathfrak{p}\big\rbrace .
\end{equation}

A. De Laire and P. Mennuni solved this minimization problem in~\cite{deLaMen1} for the nonlocal Gross-Pitaevskii equation. F. Bethuel, P. Gravejat and J.-C. Saut solved this problem for the Gross-Piteavskii equation in~\cite{BetGrSa2}. Similarly, our first result is

\begin{thm}\label{theoreme existence de travelling wave}
    Let us assume that $f\in\mathcal{C}^3(\R)$. Suppose also that $f$ satisfies the following conditions.

\begin{itemize}
    \item For all $\rho\in\R_+$,
\begin{equation}\label{hypothèse de croissance sur F minorant intermediaire}\tag{H1}
    \dfrac{c_s^2}{4} (1-\rho)^2 \leq F(\rho).
\end{equation}
    \item There exists $M\geq 0$ and $q\in [2, +\ii)$ such that for all $\rho \geq 2$,
\begin{equation}\label{hypothèse de croissance sur F majorant}\tag{H2}
    F(\rho)\leq M|1-\rho|^q.
\end{equation}
    \item \begin{equation}\label{condition suffisante pour la stabilité orbitale sur f''(1)+6f'(1)>0}\tag{H3}
    f''(1)+3f'(1)\neq 0.
\end{equation}
\end{itemize}

Then there exists $\gq_* \geq\frac{1}{32}$ such that for $\gp$ satisfying one of the following hypothesis
\begin{equation}\label{hypothèse sur gp en fonction de q_*}\tag{$H_{\gq_*}$}
    \left\{
\begin{array}{l}
    \quad\gp\in (0,\gq_*)\\
    \quad\quad\text{ or }\quad\\
    \gp=\gq_*\notin \frac{\pi}{2}+\pi\Z, \\ 
\end{array}
\right.
\end{equation}

there exists a travelling wave $\gv_c$ of speed $c\in (0,c_s)$ and of momentum $p(\gv_c)=\gp$.
\end{thm}

\begin{rem}
Note that the hypothesis \eqref{hypothèse de croissance sur F minorant intermediaire} can be reformulated. It can be stated as the existence of a positive constant $\lambda$ such that for any $\rho\in\R_+$,

\begin{equation}\label{1ere hypothèse de croissance sur F minorant intermediaire}\tag{H1'}
    \lambda (1-\rho)^2\leq F(\rho)
\end{equation}

and the fact that $c_s\leq 2\sqrt{\lambda}$. Indeed, bearing in mind that $F(1)=F'(1)=0$, we can write a Taylor expansion near 1 and observe that

\begin{equation*}
    \lambda \leq \dfrac{c_s^2}{4}\quad\text{i.e.}\quad 2 \sqrt{\lambda}\leq c_s .
\end{equation*}

This explains why, in Theorem~\ref{theoreme existence de travelling wave}, the strongest assumption \eqref{hypothèse de croissance sur F minorant intermediaire} is stated as it is. In many proofs, we shall only use the weaker assumption \eqref{1ere hypothèse de croissance sur F minorant intermediaire} when it is sufficient.
\end{rem}

\begin{rem}\label{remarque: conséquence de H1 et H2, estimée des normes}
As a consequence of~\eqref{1ere hypothèse de croissance sur F minorant intermediaire} and~\eqref{hypothèse de croissance sur F majorant}, we observe that 

$$\lambda \normLdeux{1-|u|^2}^2\leq \Vert F(|u|^2)\Vert_{L^1}\leq C' \big( 1 + \normLii{u}^{2(q-2)}\big)\Vert 1-|u|^2 \Vert_{L^2}^2.$$

Using the Sololev embedding $H^1(\R)\hookrightarrow L^\ii(\R)$, we conclude that the classical Gross-Pitaevskii energy space 

\begin{equation}\label{inclusion d'ensemble widetilde x^k dans x^k dans z^k etc}
    \big\lbrace u\in H^1_{\loc}(\R) \big| u' \in L^2 (\R), 1-|u|^2 \in L^2(\R)\big\rbrace 
\end{equation}

is exactly equal to $\energyset$. Note also that every function in $\mathcal{X}^1(\R)$ is uniformly continuous and with modulus tending to $1$. In particular, we can replace $H^1_{\loc}(\R)$ by $L^\ii(\R)$ in the definition of the energy space, and we can also take the set $\big\lbrace u\in L^{\ii}(\R)\big| u' \in L^2(\R), 1-|u| \in L^2(\R)\big\rbrace$, instead of $\energyset$. This property is specific to the one space dimension and it is not true in higher dimensions.
\end{rem}

\begin{rem}

The bound in~\eqref{theoreme de gallo condition de croissance sur f''} is a sufficient condition for~\eqref{hypothèse de croissance sur F majorant} to hold. Indeed, integrating~\eqref{theoreme de gallo condition de croissance sur f''} three times for $\rho$ large enough, we get, if $\alpha_1\notin\{1,2\}$,

\begin{equation*}
    F(\rho)\lesssim \rho^{\alpha_1}\lesssim |\rho - 1 |^{q},
\end{equation*}

with $q=\min(\alpha_1,2)$, and where $\lesssim$ means that the inequalities hold, up to a constant independent of $\rho$. If $\alpha_1=1$ or $2$, we obtain the same estimate for $\rho\geq 2$, with $q=2$ (resp. $q=3$).
\end{rem}

\begin{rem}
The quantity $f''(1)+3f'(1)$ in hypothesis~\eqref{condition suffisante pour la stabilité orbitale sur f''(1)+6f'(1)>0} is related to the constant $\Gamma$, which appears in M. Maris' article~\cite{Maris6} and in D. Chiron's (see~\cite{Chiron7,Chiron8,Chiron9}). When this number is equal to zero, the problem under consideration is known to be degenerate. The (KdV) transonic regime turns out to be a linear dispersive equation, and consequently owns no soliton.\\
\end{rem}

\begin{rem}
Regarding the above theorems, if we suppose that~\eqref{hypothèse de croissance sur F minorant intermediaire} holds true, then the zero $\xi_c$ of Theorem~\ref{existence et unicité des travelling wave si condition du zéro} (if it exists) necessarily lies in $(0,1)$. Indeed, we write that $\mathcal{N}_c(\xi)\leq \xi^2\big(c^2 - c_s^2(1-\xi)\big)$,
by~\eqref{hypothèse de croissance sur F minorant intermediaire}. So that, $0=\mathcal{N}_c(\xi_c)\leq \xi_c^2\big(c^2 -c_s^2(1-\xi_c)\big)$ and then \begin{equation*}
    \xi_c\geq 1-\dfrac{c^2}{c_s^2}\geq 0.
\end{equation*}
\end{rem}

\begin{rem}\label{remarque: exemple de non linearité f(x)=1-x+ a(1-x) ^2p-1}
    We can find explicit examples of nonlinearities $f$ that satisfy the assumptions of Theorem~\ref{existence et unicité des travelling wave si condition du zéro}, so that there exists a unique travelling wave of speed\footnote{When $c=0$, the existence and uniqueness are automatically satisfied whenever~\eqref{hypothèse de croissance sur F minorant intermediaire} holds (see Proposition~\ref{prop: la solution de tw0 est unique si int_0^1 f neq 0}) and the case $c=c_s$ cannot occur if~\eqref{condition suffisante pour la stabilité orbitale sur f''(1)+6f'(1)>0} holds (see Remark~\ref{remarque non existence de tw supersonique sonique Maris}).} $c\in (0,c_s)$, and also assumptions~\eqref{hypothèse de croissance sur F minorant intermediaire},\eqref{hypothèse de croissance sur F majorant} and~\eqref{condition suffisante pour la stabilité orbitale sur f''(1)+6f'(1)>0}, so that we can investigate the existence of minimizers for the energy and their orbital stability. Moreover, we will see in Section~\ref{section: numerical simulations} how these nonlinearities can provide quite different behaviours.

For integers $p\geq 2$, consider the function $f(\rho)=1-\rho + a(1-\rho)^{2p-1}$. We compute $c_s=\sqrt{2}$, $$F(\rho)=\dfrac{(1-\rho)^2}{2}+a_p(1-\rho)^{2p}\quad\text{with }a_p:=\dfrac{a}{2p}>0,$$ 
so that \eqref{hypothèse de croissance sur F minorant intermediaire},\eqref{hypothèse de croissance sur F majorant} and~\eqref{condition suffisante pour la stabilité orbitale sur f''(1)+6f'(1)>0} are satisfied. Now let us verify the assumptions of Theorem~\ref{existence et unicité des travelling wave si condition du zéro}. We have $\mathcal{N}_c(\xi)=\xi^2 P_p(\xi)$ with $P_p(\xi)=4a_p\xi^{2p-1}-4a_p\xi^{2p-2}+2\xi-\varepsilon^2$ and $\varepsilon^2=c_s^2- c^2$. We compute $P_p(0)=-\varepsilon^2 <0$ and $P_p(1)=c^2>0$. By the intermediate value theorem, there exists a zero $\xi_c\in (0,1)$ of $P_p$ and it is sufficient to prove that this zero is not a double root.

Indeed, if there is no double root, we can choose $\xi_c$ as the minimal root in $(0,1)$. Then we necessarily have $P'_p(\xi_c)>0$ and $P_p<0$ on $(0,\xi_c)$. In view of $\mathcal{N}_c'(\xi)=2\xi P_p(\xi)+\xi^2P'_p(\xi)$, we conclude that the same properties hold for $\mathcal{N}_c$. Let us now check that the root is single.\\
By contradiction suppose that $\xi_c$ is not single. In this case, we have $P'_p(\xi_c)=0$. Considering the variations of the polynomial function associated with $P'_p$, leads to, for all $\xi\in\R$,

\begin{equation*}
    P'_p\Big(\dfrac{2p-3}{2p-1}\Big)\leq P'_p(\xi).
\end{equation*}

In particular, for $\xi=\xi_c$, we obtain 
\begin{equation}\label{d_ap - c leq 0}
    \dfrac{1}{2a_p}-\dfrac{(2p-3)^{2p-3}}{(2p-1)^{2p-3}}\leq 0 .
\end{equation}

The nonlinearity $f$ happens to be a suitable candidate whenever we take $a<p\big(\frac{2p-1}{2p-3}\big)^{2p-3}$.
\end{rem}

Like in~\cite{deLaMen1}, the previous result of existence does not state the uniqueness of the travelling wave with a fixed momentum. The uniqueness for such general nonlinearities is difficult to establish and this question goes beyond the scope of this article. One sufficient condition to obtain the uniqueness would be a one-to-one correspondence between the speed $c$ and the momentum of the travelling wave of speed $p(\gv_c)$. Set

\begin{equation}\label{définition de S_p}
    \mathcal{S}_\gp:=\big\lbrace v\in \mathcal{N}\energyset\big| E(v)=\Emin (\gp)\text{ and }p(v)=\mathfrak{p}\big\rbrace .
\end{equation}

Theorem~\ref{theoreme existence de travelling wave} guarantees that for some $\gp$, there exists a travelling wave minimizing the energy when the momentum is fixed at $\gp$, so that $\mathcal{S}_\gp$ is not empty. In addition, we will prove that such a set is orbitally stable in the sense that we recall now.

\begin{de}\label{defintion de la stabilité orbitale}
We say that a subset $\mathcal{H}\subset \energyset$ is orbitally stable for a distance $d$ if, for any $\varepsilon >0$, there exists $\delta >0$ such that for all $\Psi_0 \in \mathcal{H}$, the solution $\Psi(t)$ of \eqref{NLS} with initial condition $\Psi_0$ satisfies the following property: if 

\begin{equation*}
    d(\Psi_0,\mathcal{H})\leq\delta,
\end{equation*}

then for all $t\in\R$, there exist $a(t),\theta(t)\in\R$ such that

\begin{equation*}
    d\Big(e^{i\theta(t)}\Psi\big(t,.-a(t)\big),\mathcal{H}\Big)\leq\varepsilon .
\end{equation*}
\end{de}

We can endow the sets in \eqref{inclusion d'ensemble widetilde x^k dans x^k dans z^k etc}, with the distance 

\begin{equation*}
    d_A(u_1,u_2):=\Vert u_1-u_2\Vert_{L^{\ii}([-A,A])} +\Vert u_1'-u_2'\Vert_{L^2}+\Vert |u_1|^2-|u_2|^2\Vert_{L^2},
\end{equation*}

where $A>0$. The corresponding metric structure is independent of the choice of the number $A$. The next theorem states the orbital stability of the minimizers of the energy when the momentum is fixed in a certain range.

\begin{thm}\label{théorème principal}
Take a function $f$ as in Theorems~\ref{global well-posedness of cauchy problem} and \ref{theoreme existence de travelling wave}. Suppose that \eqref{hypothèse de croissance sur F minorant intermediaire},\eqref{hypothèse de croissance sur F majorant} and \eqref{condition suffisante pour la stabilité orbitale sur f''(1)+6f'(1)>0} hold and choose $\gp$ satisfying \eqref{hypothèse sur gp en fonction de q_*}. Then $(\mathcal{S}_\gp,d_A)$ is orbitally stable.
\end{thm}

The first result of orbital stability for the~\eqref{NLS} equation is due to Z. Lin~\cite{LinZhiw1}, who proved it in a hydrodynamical framework given by $\big(\mathcal{NX}^1(\R),d_{\mathrm{hy}}\big)$ with

\begin{equation*}
    d_{\mathrm{hy}}(u_1,u_2)=\Vert \rho_1 - \rho_2 \Vert_{H^1} + \Vert \varphi_1'-\varphi_2'\Vert_{L^2}  + \Big|\mathrm{arg}\Big(\dfrac{u_1(0)}{u_2(0)}\Big)\Big|,\quad \quad u_j = \rho_j e^{i\varphi_j}.
\end{equation*}

It was extended by D. Chiron in~\cite{Chiron8} to $(\mathcal{X}^1(\R),d_A)$. Both results rely on a condition due to M. Grillakis, J. Shatah and W.A. Strauss~\cite{GriShSt1}, who studied the orbital stability (and instability) of solitary waves in a framework which covers, to a large extent, the one in this paper. For $c_*\in (0,c_s)$, the condition for a travelling wave of speed $c_*$ to be orbitally stable is

\begin{equation}\label{dérivée du moment en v_c par rapport à c en c=c_* est strictement négatif}
    \dfrac{dp(\mathfrak{v}_c)}{dc}_{\big|c=c_*} <0.
\end{equation}

This inequality is related to the strict concavity of the minimization curve near $c_*$. Indeed, in view of the integral expressions given for $E(\gv_c)$ and $p(\gv_c)$ in~\cite{Chiron7}, we deduce  the Hamilton group relation

\begin{align}\label{relation de groupe hamiltonien dE/dc(v_c)=cdp/dc(v_c)}
    \dfrac{dE(\mathfrak{v}_c)}{dc}=c\dfrac{dp(\mathfrak{v}_c)}{dc}.
\end{align}

Assuming that~\eqref{dérivée du moment en v_c par rapport à c en c=c_* est strictement négatif} holds for one $c_*\in(0,c_s)$ and using the inverse function theorem, we obtain the expression of the energy in terms of the momentum
\begin{equation}\label{le signe de mathcal E'' est le meme que celui de p'}
    \dfrac{d^2\mathcal{E}}{dp^2}\big(p(\mathfrak{v}_c)\big)=\Big(\dfrac{dp(\mathfrak{v}_c)}{dc}\Big)^{-1},
\end{equation}
where we have set $\mathcal{E}\big(p(\gv_c)\big)=E(\gv_c)$. This relates the sign of $\frac{dp(\mathfrak{v}_c)}{dc}$ and the concavity of the function $\mathcal{E}$, that we will establish in the sequel.

The assumptions of Theorem~\ref{théorème principal} do not mention any condition like~\eqref{dérivée du moment en v_c par rapport à c en c=c_* est strictement négatif}, unlike it was the case in~\cite{Chiron7,GriShSt1}. Our hypothesis are more suitable than \eqref{dérivée du moment en v_c par rapport à c en c=c_* est strictement négatif} in the sense that the class of functions $f$ for which we have orbital stability is more explicit. Here we only make elementary assumptions on $f$ and we adapt the variational method in~\cite{BetGrSa2} in order to prove the orbital stability of $\mathcal{S}_\gp$. Unlike Theorem~\ref{existence et unicité des travelling wave si condition du zéro}, this variational approach is also expected to give existence results for higher dimensions (see~\cite{BetGrSa1,BetOrSm1}). Here we show the existence of a branch of stable travelling waves while weaning off both the assumptions of Theorem~\ref{existence et unicité des travelling wave si condition du zéro} and the condition of Grillakis, Shatah and Strauss~\eqref{dérivée du moment en v_c par rapport à c en c=c_* est strictement négatif}. Moreover, we have the explicit lower bound $\gq_*\geq\frac{1}{32}$ on the length of the branch of stable solitons.

More precisely, we prove at the same time the existence and the orbital stability of minimizers for the energy when the momentum is fixed in a certain range~\eqref{hypothèse sur gp en fonction de q_*}. However, when we fix the speed $c\in (0,c_s)$, we cannot prove the orbital stability of the travelling wave associated with this speed $c$. This is due to the fact that we cannot prove the uniqueness of the travelling wave minimizing $\Emin$ at fixed momentum. Because of this fact, the concentration-compactness argument yields a travelling wave whose speed is not explicitly given with respect to the constraint $\gp$. Up to more restrictive assumptions for the nonlinearity $f$, we would expect that the travelling waves are orbitally stable for a certain range of speed $c\in(c_{*},c_s)$. And this, by showing that $c\mapsto p(\gv_c)$ is smooth enough and that \eqref{dérivée du moment en v_c par rapport à c en c=c_* est strictement négatif} holds for speeds close to $c_s$.

The proof of Theorems~\ref{theoreme existence de travelling wave} and~\ref{théorème principal} relies on the variational method introduced for the first time in~\cite{CazeLio1} and then applied in many articles~\cite{BetGrSa2,BeGrSaS1,ChirMar2,deLaMen1}.
It is based on a concentration-compactness argument for the study of the minimization of the energy when the momentum is fixed.\\
The minimizing energy is proved to be concave and strictly sub-additive, which ultimately provides compactness, and then orbital stability.
Throughout this paper, we shall suppose that the function $f$ is at least $\mathcal{C}^3(\R)$ and that $f'''$ is bounded. These assumptions are crucial to obtain some of the following estimates (for instance the right-hand side of \eqref{prop: inegalité a droite et a gauche de Emin par des puissances de q}) which are required by the concentration-compactness argument of Theorem~\ref{theoreme de continuité de l'energie dans le cas d'une non linéarité générale}.

\subsection{Sketch of the proofs}
We now give the main steps of the argument. Especially, we explain how the hypothesis of Theorem~\ref{theoreme existence de travelling wave}, \eqref{hypothèse de croissance sur F minorant intermediaire}, \eqref{hypothèse de croissance sur F majorant} and \eqref{condition suffisante pour la stabilité orbitale sur f''(1)+6f'(1)>0}, yield the existence of a family of travelling waves and its orbital stability when it is parametrized by a momentum satisfying \eqref{hypothèse sur gp en fonction de q_*}. The main theorem relies on the variational interpretation of \eqref{TWC} as a minimization problem under constraints, for which we can establish the compactness of the minimizing sequences by a concentration-compactness argument. Knowing those facts, we prove the orbital stability (à la Cazenave-Lions) by contradiction. Assuming that the set $\mathcal{S}_\gp$ is not orbitally stable, we exhibit a pseudo-minimizing sequence for the variational problem that tends to $\gv_c$ and this brings the contradiction.\\

Above all, we begin by defining properly the momentum $p$. This quantity is known to have a rigorous sense (see~\cite{BetGrSa2,BeGrSaS1}) on the set $\Nenergyset$. When a function $u$ in this set is lifted as $u=\rho e^{i\varphi}$, it is given by the formula

\begin{equation*}
    p(u)=\dfrac{1}{2}\int_\R (1-\rho^2)\varphi '.
\end{equation*}

We are then allowed to study in details the properties of $\Emin$ whose behaviour is similar to the one of the minimizing energy in the Gross-Pitaevskii case. One of the first property to notice is the fact that $\Emin$ is even. Accordingly to this, we display the graph of $\Emin$ on $\R_+$ in the Gross-Pitaevskii case.

\begin{center}
    \includegraphics[scale=0.4]{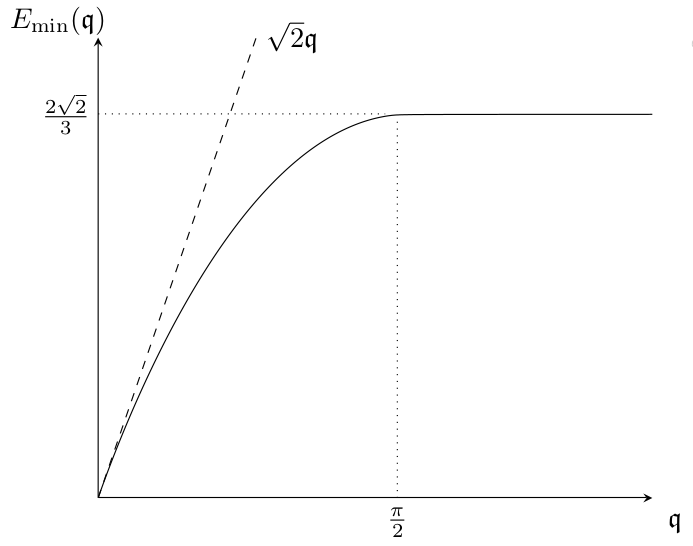}\\
    \small Graph of $E_{\min}$ for the Gross-Pitaevskii nonlinearity $f(\rho)=1-\rho$ (graph from~\cite{deLaMen1}).
\end{center}

We observe that a transition occurs at $\frac{\pi}{2}$, which we aim at understanding. We introduce the quantity 

\begin{equation*}
    \gq_* = \sup\big\lbrace \gq >0\big| \forall v\in \energyset, E(v)\leq \Emin(\gq)\Rightarrow \inf_\R |v| >0\big\rbrace .
\end{equation*}

\begin{prop}\label{prop Emin est lipschitz continuous, inégalité de lipschitz, inegalité a droite et a gauche de Emin par des puissances de q}
Let us assume that \eqref{hypothèse de croissance sur F minorant intermediaire} and \eqref{hypothèse de croissance sur F majorant} are satisfied. Then the following statements hold.

    $(i)$ The function $E_{\min}$ is nonnegative, even, Lipschitz continuous on $\R$, with
        
    \begin{equation*}\label{prop: inegalité de type lipschitz sur Emin}
        |E_{\min}(\mathfrak{p})-E_{\min}(\mathfrak{q})|\leq c_s |\mathfrak{p}-\mathfrak{q}|,\text{ for all }\mathfrak{p},\mathfrak{q}\in\R.
    \end{equation*}

    $(ii)$ $\Emin$ is nondecreasing and concave on  $\R_+$ and strictly increasing on $[0,\gq_*]$.

    $(iii)$ Suppose that \eqref{condition suffisante pour la stabilité orbitale sur f''(1)+6f'(1)>0} holds. There exist positive constants $\mathfrak{q}_0, K_0,K_1,K_2$ such that

\begin{equation}\label{prop: inegalité a droite et a gauche de Emin par des puissances de q}
    c_s\mathfrak{q}-K_0 \mathfrak{q}^{\frac{5}{3}}\leq E_{\min}(\mathfrak{q})\leq c_s\mathfrak{q}-K_1\mathfrak{q}^{\frac{5}{3}}+K_2\mathfrak{q}^\frac{7}{3}\text{ for all }\mathfrak{q}\in [0,\mathfrak{q}_0].
    \end{equation}

    $(iv)$ Suppose that \eqref{condition suffisante pour la stabilité orbitale sur f''(1)+6f'(1)>0} holds. Then $\Emin$ is strictly subadditive.
\end{prop}

\begin{rem}
Here concavity means that the function $\Emin$ satisfies the property of concavity in a large sense. When it is relevant, we will distinguish concave and strictly concave.
\end{rem}

The proof follows from constructing a sequence of test functions which approximates the infimum $\Emin(\gq)$ for any $\gq\in\R$. We obtain for instance $(i)$ by working on such sequences and letting $n$ tend to $+\ii$. The most significant property is the inequality in the right-hand side of \eqref{prop: inegalité a droite et a gauche de Emin par des puissances de q}. For deriving it, we construct a family of functions in the transonic regime $c\rightarrow c_s$, the behaviour of which resembles the Korteweg-De Vries solutions due to condition~\eqref{condition suffisante pour la stabilité orbitale sur f''(1)+6f'(1)>0}.

We shall put a special interest in the travelling wave of speed $c=0$. Indeed, we shall see that this travelling wave is the only one that vanishes on $\R$. Especially, it is a minimizer for the energy among the functions that vanishes and this will help us understand why some values are prescribed for the momentum.

The minimization problem $\Emin (\gp)$ with $\gp\in (0,\gq_*)$ is attained on solutions $\mathfrak{v}_c$ of \eqref{TWC} satisfying $p(\mathfrak{v}_c)=\gp$. This is also true whenever $\gp=\gq_*$ with $\gq_* \neq \frac{\pi}{2}\mod \pi$. To check this claim, the properties of $\Emin$ can be used to prove that we can find a minimizing sequence that converges (in a sense to be precised and up to a subsequence) to $\mathfrak{v}_c$. This compactness result follows from a concentration-compactness argument (see~\cite{Lions1}) and is stated in our framework as follows.

\begin{thm}\label{theoreme de continuité de l'energie dans le cas d'une non linéarité générale}
Assume that \eqref{hypothèse de croissance sur F minorant intermediaire}, \eqref{hypothèse de croissance sur F majorant}, \eqref{condition suffisante pour la stabilité orbitale sur f''(1)+6f'(1)>0} and \eqref{hypothèse sur gp en fonction de q_*} hold. Let $(u_n)$ be a pseudo-minimizing sequence i.e. satisfying 

\begin{equation*}\label{etre une suite pseudo minimisante}
    p(u_n)\ntend\mathfrak{p}\quad\text{and}\quad E(u_n)\ntend \Emin (\mathfrak{p}).
\end{equation*}

Then, there exist a subsequence $(u_{\sigma(n)})$, a sequence of points $(a_n)$, a real number $\theta$ and a non constant solution $\mathfrak{v}_c$ of \eqref{TWC} such that

\begin{equation*}
    u_{\sigma(n)}(.+a_{\sigma(n)})\ntend e^{i\theta} \mathfrak{v}_{c}\quad\text{in }\mathcal{C}^0_{\loc}(\R).
\end{equation*}

Moreover, we have

\begin{equation*}
    u'_{\sigma(n)}(.+a_{\sigma(n)})\ntend e^{i\theta} \mathfrak{v}'_{c}\quad\text{in }L^2(\R),
\end{equation*}

and

\begin{equation*}
    F\big(|u_{\sigma(n)}(.+a_{\sigma(n)})|^2\big) \ntend   F\big(|\mathfrak{v}_{c}|^2\big)   \quad   \text{in }L^1(\R).
\end{equation*}

In addition, \begin{equation}\label{gv_c in S_gp}
\gv_c\in\mathcal{S}_\gp.
\end{equation}
\end{thm}

This theorem is a consequence of a general phenomenon that occurs in minimization problems with constraints and was highlighted by P.-L. Lions~\cite{Lions1}. The concentration-compactness theorem states that the minimization sequences are compact (up to the invariances) if vanishing is forbidden and some sub-additivity inequality is strict. The quantity involved in the sub-additivity is the infimum of the problem considered as a function of the value of the constraint. In our framework, it is the infimum of the energy when the momentum is fixed. In other words, the strict sub-additive property of $\Emin$ (Proposition~\ref{prop Emin est lipschitz continuous, inégalité de lipschitz, inegalité a droite et a gauche de Emin par des puissances de q}) is a sufficient condition for Theorem~\ref{theoreme de continuité de l'energie dans le cas d'une non linéarité générale} to hold.

Section~\ref{section propriétés de la courbe minimisante d'énergie} is devoted to the properties of the minimization curve $\Emin$. Section~\ref{section charcteristics of the kink solution} deals with the properties of the kink solution for $c=0$. In Section \ref{section méthode variationnelle, minimisation de l'énergie à moment fixé}, we see how the variational method under the constraint that the momentum is fixed can be implemented to a general nonlinearity $f$. Section~\ref{section stabilité orbitale} deals with the orbital stability of the minimizers of this variational problem. Finally, in Section~\ref{section: numerical simulations}, we display several numerical simulations, that lay emphasis on the different behaviours that can occur according to the nonlinearities. In particular, we give examples where~\eqref{hypothèse de croissance sur F minorant intermediaire},\eqref{hypothèse de croissance sur F majorant} and~\eqref{condition suffisante pour la stabilité orbitale sur f''(1)+6f'(1)>0} are satisfied with $\gq_*=\frac{\pi}{2}$ as in the Gross-Pitaevskii case, but also with $\gq_*<\frac{\pi}{2}$ and $\gq_* >\frac{\pi}{2}$.

\section{Properties of the minimization curve}\label{section propriétés de la courbe minimisante d'énergie}
For the study of the minimization curve, we introduce the set

\begin{equation*}
    \mathcal{NX}^\ii_0(\R)=\big\lbrace v\in\mathcal{N}\energyset\cap \mathcal{C}^\ii(\R)\big| \exists R >0\text{ s.t. }v\text{ is constant on }(-R,R)^c \big\rbrace .
\end{equation*}

The next result shows that $\Emin$ is well defined and its graph lies under the line $E=c_s p$ on $\R_+$.

\begin{lem}\label{lemme qui definit bien Emin et qui donne une borne par c_s q}
For all $\mathfrak{q}\in\R$, there exists a sequence $(v_n)\in \mathcal{NX}^\ii_0(\R)^\N$ satisfying
\begin{equation}\label{lemme: il existe une suite de moment q dont l'energie tend vers c_s q}
    p(v_n)=\mathfrak{q}\quad\text{and}\quad E(v_n)\ntend c_s|\mathfrak{q}|.
\end{equation}

In particular the function $E_{\min}:\R\rightarrow\R_+$ is well-defined, and for all $\mathfrak{q}\geq 0$, 

\begin{equation}\label{lem: Eminq leq cs q}
    0\leq E_{\min}(\mathfrak{q})\leq c_s\mathfrak{q}.
\end{equation}
\end{lem}

\begin{proof}
The case $\mathfrak{q}=0$ results from taking $v\equiv 1\in  \mathcal{NX}^\ii_0(\R)$. Now let us assume that $\mathfrak{q}> 0$ and consider $\chi\in\mathcal{C}^{\ii}_{c}(\R)$ such that $\int_\R (\chi')^2=\frac{\mathfrak{q}}{c_s}$. Let us also define

\begin{equation*}
    a:=\int_\R (\chi')^3,\quad\alpha_n:=\dfrac{1}{n}\quad\text{and}\quad\beta_n:=\dfrac{1}{n^2}\Big( 1-\dfrac{c_s}{2\mathfrak{q}n}a\Big).
\end{equation*}

Then we set \begin{equation*}
    v_n:=\rho_n e^{i\varphi_n},\quad\text{where }\rho_n(x):=1-\alpha_n\chi'(\beta_n x)\text{ and }\varphi_n(x):=c_s\dfrac{\alpha_n}{\beta_n}\chi(\beta_n x).
\end{equation*}

Because of the asymptotic properties of both sequences $(\alpha_n)$ and $(\beta_n)$, the function $v_n$ is well-defined and does not vanish for $n$ large enough. Thus its momentum is also well-defined and we have
\begin{align*}
    p(v_n) & = \dfrac{1}{2}\int_\R (1-\rho_n^2 ) \varphi_n '\\
    & = \dfrac{1}{2}\int_\R \big( 2\alpha_n \chi'(\beta_n x) - \alpha_n^2 \chi'(\beta_n x) ^2 \big)c_s \alpha_n \chi'(\beta_n x ) dx \\
    & =\dfrac{\alpha_n^2}{\beta_n}\mathfrak{q}-\dfrac{c_s \alpha_n^3}{2\beta_n}a=\mathfrak{q}.
\end{align*}

Computations similar as for the momentum show that
\begin{align*}
    E_k(v_n) & =\dfrac{1}{2}\int_\R \rho_n'^2+\rho_n^2 \varphi_n'^2\\
    &=\dfrac{\alpha_n^2}{\beta_n}\int_\R (1-\alpha_n\chi')(\chi')^2+\dfrac{\alpha_n^2\beta_n}{2}\int_\R (\chi'')^2\\
    &\ntend \dfrac{c_s^2}{2}\int_\R (\chi')^2=\dfrac{c_s \mathfrak{q}}{2}.
\end{align*}

For the potential energy, we use a Taylor expansion with integral remainder: for all $x\in\R$,
\begin{equation*}
    F(1+x)=\dfrac{c_s^2}{4}x^2-\dfrac{x^3}{2}\int_0^1 (1-t)^2 f''(1+xt)dt.
\end{equation*}

We now compute the limit of $E_p(v_n)$. Since $\chi'$ is compactly supported and $f''$ is continuous, there exists $M$ independent of $x$ and $n$ such that $\big|f''\big(1-2t\alpha_n \chi'(\beta_n x)+t\alpha_n^2 \chi'(\beta_n x)^2\big)\big|\leq M$. Therefore, replacing in the previous Taylor formula $x$ by $-2\alpha_n \chi'(\beta_n x)+\alpha_n^2 \chi'(\beta_n x)^2$, and introducing the resulting expression in the integrand, we obtain
\begin{align*}
    E_p(v_n)&=\dfrac{1}{2}\int_\R F\big(1-2\alpha_n \chi'(\beta_n x)+\alpha_n^2 \chi'(\beta_n x)\big)dx \\
    &=\dfrac{c_s^2}{8\beta_n}\int_\R (\alpha_n^2 \chi'^2-2\alpha_n\chi')^2+\petitodeunnii.
\end{align*}

The second equality in this formula follows from the fact that $\chi\in\mathcal{C}_{c}^{\ii}(\R)$ and the asymptotic properties corresponding to the sequences $(\alpha_n)$ and $(\beta_n)$. Since 

\begin{equation*}
    \dfrac{\alpha_n^2}{\beta_n}\ntend 1,\quad\text{while }\dfrac{\alpha_n^3}{\beta_n},\dfrac{\alpha_n^4}{\beta_n}\ntend 0,
\end{equation*}

we obtain

\begin{equation*}
    E_p(v_n)\ntend \dfrac{c_s^2}{8}\int_\R 4\chi'^2=\dfrac{c_s\mathfrak{q}}{2}.
\end{equation*}

Therefore we conclude that (\ref{lemme: il existe une suite de moment q dont l'energie tend vers c_s q}) holds true for $\mathfrak{q}\geq 0$. In the case $\mathfrak{q}<0$, it is enough to proceed as above taking

\begin{equation*}
    \int_\R \chi'^2 =-\dfrac{\mathfrak{q}}{c_s}\quad\text{ and }v_n=\rho_n e^{-i\varphi_n}.
\end{equation*}

\end{proof}

In view of Lemma~\ref{lemme qui definit bien Emin et qui donne une borne par c_s q}, it does not matter to define $E_{\min}$ as the minimizer on the set $\mathcal{NX}^1(\R)$ or on $\mathcal{N}\mathcal{X}^\ii_0(\R)$. Therefore $E_{\min}$ is well-defined and is moreover even due to the next lemma which is inspired of~\cite{deLaMen1}.

\begin{lem}
$\Emin$ is even.
\end{lem}

\begin{proof}
Let $\gq \in\R$ and $u_n= \rho_n e^{i\varphi_n}\in \Nenergyset$ such that $E(u_n)\ntend \Emin(\gq)$ and $p(u_n)=\gq$. We set $v_n = \rho_n e^{-i\varphi_n}$. We verify that $E(v_n)=E(u_n)\ntend \Emin(\gq)$ and $p(v_n)=-p(u_n)=-\gq$. As a consequence, we have $E(v_n)\geq \Emin(-\gq)$. Now letting $n\rightarrow + \ii$, we obtain $\Emin(\gq) \geq \Emin(-\gq)$ and the reverse inequality follows from replacing $\gq$ by $-\gq$. We conclude that $\Emin$ is even.
\end{proof}

Concerning the density of the space $\mathcal{NX}^\ii_0(\R)$ in $\Nenergyset$, we have the following result.
\begin{lem}\label{lemme densité de Eii dans NE}
Assume that \eqref{hypothèse de croissance sur F majorant} holds. Let $v=\rho e^{i\varphi}\in \Nenergyset$. Then there exists a sequence of functions $(v_n)=(\rho_n e^{i\varphi_n})\in \mathcal{NX}^\ii_0(\R)^\N$, with $\rho_n - 1,\varphi_n\in\mathcal{C}_c^{\ii}(\R)$, such that 

\begin{equation}\label{lemme: rho_n tend vers rho H1, theta_n' tend vers theta' L2}
    \Vert \rho_n -\rho\Vert_{H^1}+\Vert\varphi_n' -\varphi'\Vert_{L^2}\ntend 0.
\end{equation}

In particular, 

\begin{equation}\label{lemme convergence de E(v_n) vers E(v) et de p(v_n) vers p(v) quand convergence de rho_n dans H^1 et de theta_n' dans L^2}
    E(v_n)\ntend E(v)\quad\text{and}\quad p(v_n)\ntend p(v).
\end{equation}
\end{lem}

\begin{proof}
The existence of $(\rho_n)$ and $(\varphi_n)$ such as in the lemma and satisfying \eqref{lemme: rho_n tend vers rho H1, theta_n' tend vers theta' L2} was shown in~\cite[Lemma 3.4]{deLaMen1}. Hence we just show the convergences \eqref{lemme convergence de E(v_n) vers E(v) et de p(v_n) vers p(v) quand convergence de rho_n dans H^1 et de theta_n' dans L^2}.\\
We have
\begin{align*}
    \Vert(1-\rho_n^2)-(1-\rho^2)\Vert_{L^2}&\leq\Vert \rho_n-\rho\Vert_{L^2} \big(\Vert\rho_n\Vert_{L^{\ii}} + \Vert \rho\Vert_{L^{\ii}}\big).
\end{align*}

In view of (\ref{lemme: rho_n tend vers rho H1, theta_n' tend vers theta' L2}), and because of the Sobolev embedding $H^1(\R)\hookrightarrow L^\ii(\R)$, the norms $\Vert \rho_n\Vert_{L^{\ii}}$ are uniformly bounded with respect to $n$ and then

\begin{equation*}
    1-\rho_n^2\ntend 1-\rho^2\quad\text{in }L^2(\R).
\end{equation*}

This and the convergence of $\varphi_n'$ to $\varphi'$ in $L^2(\R)$ imply that \begin{equation*}
    p(v_n)\ntend p(v).
\end{equation*}

Secondly, we write
\begin{equation*}
    E_k(v_n)=\dfrac{1}{2}\int_\R (\rho'_n)^2+\rho_n^2 (\varphi_n')^2.
\end{equation*}

We have
\begin{equation*}
    \int_\R \rho'^2_n\ntend \int_\R \rho'^2,
\end{equation*}

and 
\begin{align*}
    \Vert \rho_n\varphi'_n-\rho\varphi'\Vert_{L^2} &\leq \Vert \rho_n-\rho\Vert_{L^{\ii}} \Vert\varphi_n'\Vert_{L^2}+\Vert \rho\Vert_{L^{\ii}} \Vert\varphi_n'-\varphi'\Vert_{L^2}.
\end{align*}

The same arguments lead to 
\begin{equation*}
    \rho_n\varphi'_n\ntend \rho\varphi'\quad\text{in } L^2(\R),
\end{equation*}

so that
\begin{equation*}
    E_k(v_n)\ntend E_k(v).
\end{equation*}

Finally we address the convergence of $E_p(v_n)$. Since \begin{equation*}
    1-\rho_n \ntend 1-\rho\quad\text{in }H^1(\R),
\end{equation*}

by Lemma~\ref{lemme convergence de 1-rho_n dans H^1 implique convergence de F(rho_n^2)} in the Appendix, we obtain the convergence of $E_p(v_n)$ to $E_p(v)$. 
\end{proof}

We can modify a function with energy close to $E_{\min}(\mathfrak{q})$ such that it is constant far away, but the momentum remains unchanged. This property implies the continuity of $E_{\min}$. We refer to~\cite[Corollary 3.7, Corollary 3.8, Proposition 3.9]{deLaMen1} for the proofs of the three following results\footnote{There are small modifications in the proofs of these results. They are based on estimating the potential energy of a piecewise smooth function. In our case and by the jump formula, the potential energy of a piecewise smooth function is the sum of the potential energies of this function on each branch. Whereas in~\cite{deLaMen1}, it is necessary to deal with a nonlocal potential and this requires additional assumptions to have the good estimate.}.

\begin{cor}\label{corollaire il existe une suite dans Eii de moment p(u) et dont l'énergie tend vers E(u)}
Let $u=\rho e^{i\varphi}\in \Nenergyset$. There exists a sequence $(u_n)\in \mathcal{NX}^\ii_0(\R)^\N$ such that 

\begin{equation*}
    p(u_n)=p(u)\quad\text{and}\quad E(u_n)\ntend E(u).
\end{equation*}
\end{cor}

\begin{cor}\label{corollaire l'inf de Emin peut se prendre sur Eii}
For all $\mathfrak{q}\geq 0$ and $\varepsilon >0$, there exists a function $v\in \mathcal{NX}^\ii_0(\R)$ such that 

\begin{equation*}
    p(v)=\mathfrak{q}\quad\text{and}\quad E(v)<E_{\min}(\mathfrak{q})+\varepsilon.
\end{equation*}
In particular, 

\begin{equation*}
    E_{\min}(\mathfrak{q})=\inf\big\lbrace E(v)\big| v\in  \mathcal{NX}^\ii_0(\R),p(v)=\mathfrak{q}\big\rbrace.
\end{equation*}

\end{cor}

Now we state that $\Emin$ is a $c_s$-Lipschitz function on $\R_+$.
\begin{prop}
$E_{\min}$ is continuous and 
\begin{equation*}
    |E_{\min}(\mathfrak{p})-E_{\min}(\mathfrak{q})|\leq c_s |\mathfrak{p}-\mathfrak{q}|\quad\text{for all }\gp,\gq\geq 0.
\end{equation*}
\end{prop}

We now address the concavity of the function $\Emin$, by appropriating the proof of Lemma~3.5 in~\cite{BetGrSa1}.

\begin{prop}\label{proposition Emin est concave}
$E_{\min}$ is a concave function on $\R_+$.
\end{prop}

\begin{proof}

Let $u=\rho e^{i\varphi}\in \mathcal{NX}^\ii_0(\R)$, such that

\begin{equation*}
    p(u)=\dfrac{1}{2}\int_\R (1-\rho^2)\varphi' =: \dfrac{\mathfrak{p}+\mathfrak{q}}{2}.
\end{equation*}

We set \begin{equation*}
    Q(a)=\dfrac{1}{2}\int_{-\ii}^a (1-\rho^2)\varphi',
\end{equation*}

\begin{equation*}
    u_-(x) =\left\{
\begin{array}{l}
    u(x)\text{ if }x\leq a \\
    e^{i\varphi_a}\overline{u(2a-x)}\text{ if }x>a
    \end{array}
\right. \quad \text{and}\quad
    u_+(x)=\left\{
\begin{array}{l}
    e^{i\varphi_a}\overline{u(2a-x)}\text{ if }x\leq a \\
    u(x)\text{ if }x > a,
\end{array}
\right.
\end{equation*}

with $\varphi_a\in\R$ such that $u_\pm$ is continuous in $a$. Without loss of generality, we can suppose that $\mathfrak{p}\leq\mathfrak{q}$. Moreover, by continuity of $Q$ on $\R$, the intermediate value theorem provides $a\in\R$ such that $Q(a)=\frac{\mathfrak{p}}{2}$. Since we have $p(u_-)=2Q(a)=\mathfrak{p}$ and $p(u_+)=2\big(p(u)-Q(a)\big)=\mathfrak{q}$, we obtain

\begin{equation*}
    \Emin(\gp)+\Emin(\gq) \leq E(v_+)+E(v_-)=2E(u),
\end{equation*}

so that \begin{equation*}
    \dfrac{E_{\min}(\mathfrak{p})+E_{\min}(\mathfrak{q})}{2}\leq E(u).
\end{equation*}

Since the choice of $u$ such that $p(u)=\frac{\mathfrak{p}+\mathfrak{q}}{2}$ is arbitrary, we obtain

\begin{equation*}
    \dfrac{E_{\min}(\mathfrak{p})+E_{\min}(\mathfrak{q})}{2}\leq E_{\min}\Big(\dfrac{\mathfrak{p}+\mathfrak{q}}{2}\Big),
\end{equation*}

which shows, by continuity of $\Emin$, that $E_{\min}$ is concave.

\end{proof}

Regarding the monotonicity of $\Emin$, we state 

\begin{prop}\label{lemme Emin est croissante}
$\Emin$ is nondecreasing on $\R_+$.
\end{prop}

\begin{proof}
We take $0 <\gp <\gq$ and $\lambda :=\frac{\gp}{\gq}\in (0,1)$. For $\delta >0$, we take $v=\rho e^{i\varphi}\in\Nenergyset$ such that $E(v) < \Emin (\gq) +\delta$ and $p(v)=\gq$. Then the function $v_\lambda = \rho e^{i\lambda\varphi}$ satisfies $p(v_\lambda)=\lambda \gq=\gp$ and $E(v_\lambda)\leq E(v)$. Therefore 

\begin{equation*}
    \Emin (\gp)\leq E(v_\lambda)\leq E(v) < \Emin(\gq) +\delta .
\end{equation*}

The conclusion follows letting $\delta\rightarrow 0$.\\
\end{proof}

Recall that $\gq_*=\sup\{\gq>0 | \forall v\in \energyset, E(v)\leq\Emin(\gq)\Rightarrow \inf_\R |v| >0\}$.
\begin{prop}\label{lemme Emin est strictement croissante}
$\Emin$ is strictly increasing on $[0,\gq_*]$.
\end{prop}

The proof relies on a special property of the black soliton $\gv_0$. Therefore we will prove this in Section~\ref{section charcteristics of the kink solution}. Nonetheless, we can give a lower bound for $\gq_*$.

\begin{prop}
    Assume that~\eqref{hypothèse de croissance sur F minorant intermediaire} holds. We have 
    \begin{equation*}
        \gq_*\geq \dfrac{1}{32}.
    \end{equation*}
\end{prop}

\begin{proof}
    Set $\gq >0$ and $v\in\energyset$ such that $E(v)\leq \Emin(\gq)$. Then by estimate~\eqref{preuve lemme appendix normlii eta inf (a+sqrt a(a+4))/2} in the appendix, combined with the estimate~\eqref{lem: Eminq leq cs q}, and the fact that $c_s=2\sqrt{\lambda}$, we obtain
    $$\normLii{1-|v|^2}\leq 8\gq + 4\sqrt{4\gq^2+\gq}.$$

Labelling the right-hand side of this estimate $h(\gq)=8\gq + 4\sqrt{4\gq^2+\gq}$, we obtain a strictly increasing function on $\R_+$. There exists a unique real number $\gq$ satisfying $h(\gq)=1$, which is larger than $\frac{1}{32}$. Thus, by definition of $\gq_*$, it implies that $\gq_*\geq\frac{1}{32}.$
\end{proof}

Now, we give estimates on the minimization curve near $0$ which are crucial to prove the strict subadditive property of $\Emin$.

\begin{prop}\label{prop inegalité a gauche de Emin par des puissances de q}
Assume that~\eqref{1ere hypothèse de croissance sur F minorant intermediaire} holds. Then there exist constants $\mathfrak{q}_0,K_0 >0$ such that 
\begin{equation*}
    c_s \mathfrak{q} -K_0\mathfrak{q}^{\frac{5}{3}}\leq \Emin(\gq),\quad\text{for all }\mathfrak{q}\in [0,\mathfrak{q}_0].
\end{equation*}
\end{prop}

\begin{proof}
Invoking Corollary \ref{corollaire l'inf de Emin peut se prendre sur Eii} and Lemma~\ref{lemme qui definit bien Emin et qui donne une borne par c_s q}, for $\delta \in (0,\frac{1}{2})$, there exists $v\in \Nenergyset$ such that $p(v)=\mathfrak{q}$ and $E(v)  < \Emin(\gq) +\delta\leq c_s \mathfrak{q}+\delta$. Then, using Lemma~\ref{lemme estimée de la norme infinie de eta en fonction de l'energie E(v)} in the appendix, we conclude that there exist $\mathfrak{q}_0 >0$ small enough and a constant $K >0$ such that, for $\gq\leq \gq_0$, $E(v)\leq 1$ and 

\begin{equation*}
    \big| 1-|v|^2\big|\leq K(c_s \gq +\delta).
\end{equation*}

We can assume that $K(c_s \gq +\delta) <1$, therefore we can apply inequality \eqref{cor: 2eme inegalité} in Corollary~\ref{corollaire du lemme sqrt(f rho^2)|varphi'| leq lambda e(v)/rho} to obtain

\begin{equation*}
    \sqrt{4\lambda(1-K(c_s \gq +\delta))}|p(v)|\leq E(v) <\Emin(\gq) +\delta,
\end{equation*}

i.e.

\begin{equation*}
    c_s  \gq \sqrt{1-K(c_s\gq +\delta)}< \Emin(\gq)+\delta.
\end{equation*}

Up to taking smaller numbers $\delta,\mathfrak{q}_0$, we can suppose that 
\begin{equation*}
    1-\varepsilon^{\frac{2}{3}}\leq\sqrt{1-\varepsilon},
\end{equation*}
with $\varepsilon = K(c_s \gq+\delta)$.
We finally infer that there exists $K_0 >0$ such that
\begin{equation*}
    c_s\mathfrak{q}-K_0\mathfrak{q}^{\frac{5}{3}}\leq \Emin(\gq),\quad\text{for all }\mathfrak{q}\leq\mathfrak{q}_0.
\end{equation*}
\end{proof}

\begin{prop}\label{prop inegalité a droite de Emin par des puissances de q}
Assume that \eqref{condition suffisante pour la stabilité orbitale sur f''(1)+6f'(1)>0} holds. There exist constants $\mathfrak{q}_1,K_1,K_2 >0$ depending on $f$, such that

\begin{equation}\label{prop: inegalité a droite de Emin par des puissances de q}
    \Emin(\gq) \leq c_s\mathfrak{q}-K_1\mathfrak{q}^{\frac{5}{3}}+K_2\mathfrak{q}^{\frac{7}{3}},\quad\text{for all }\gq\in [0,\gq_1].
\end{equation}

\end{prop}

To obtain this estimate, we introduce special test functions $v_\varepsilon$ for which this inequality holds in the limit $\varepsilon\rightarrow 0$. For that purpose, we use the ansatz

\begin{equation*}
    u_\varepsilon=\big(1+\varepsilon^2 A_\varepsilon(\varepsilon .)\big)e^{i\varepsilon\varphi_\varepsilon(\varepsilon.)},
\end{equation*}

where $A_\varepsilon,\varphi_\varepsilon$ are supposed to be smooth, with bounded derivatives. This ansatz comes from the fact that the Korteweg-De Vries equation provides a good approximation of the Gross-Pitaevskii equation in the long wave regime $\varepsilon\rightarrow 0$ (see~\cite{KuznZak1,BeGrSaS2,BeGrSaS3,ChirRou2,Chiron9} for more details on this point). First we perform some formal computations, assuming that this ansatz is a solution in order to find the better choice for the functions  $A_\varepsilon$ and $\varphi_\varepsilon$. Then we compute its momentum and energy and we complete the proof of the proposition.

Assuming that $u_\varepsilon$ is a solution of \eqref{TWC} leads us to the system

\begin{equation}\label{premier systeme TWc en remplacant u par l'ansatz u_eps}
    \left\{
\begin{array}{l}
    2\varepsilon^2 A_\varepsilon'\varphi_\varepsilon'+(1+\varepsilon^2 A_\varepsilon)\varphi_\varepsilon''+cA_\varepsilon'=0, \\
    \varepsilon^4 A_\varepsilon''-\varepsilon^4 (\varphi_\varepsilon')^2(1+\varepsilon^2 A_\varepsilon)-c\varepsilon^2 (1+\varepsilon^2 A_\varepsilon)\varphi_\varepsilon'+(1+\varepsilon^2 A_\varepsilon)f\big( (1+\varepsilon^2 A_\varepsilon)^2\big)=0.
    \end{array}
\right.
\end{equation}

We have

\begin{equation*}
    f\big( (1+\varepsilon^2 A_\varepsilon)^2\big)=2f'(1)A_\varepsilon\varepsilon^2 + \big( f'(1)+2f''(1)\big)A_\varepsilon^2  \varepsilon^4+ R_\varepsilon(A_\varepsilon)\varepsilon^6,
\end{equation*}

with 
\begin{equation*}
    | R_\varepsilon(z)| \leq C\Vert f'''\Vert_{L^{\ii}} |z|^3,
\end{equation*}

provided that $f'''$ is bounded. As a matter of fact, we are going to see that $f'''$  need not be bounded on the whole line but only locally bounded, which is the case since $f\in\mathcal{C}^3(\R)$. In this way, we formally deduce from (\ref{premier systeme TWc en remplacant u par l'ansatz u_eps}) that

\begin{equation}\label{deuxieme systeme TWc en remplacant u par l'ansatz u_eps}
    \left\{
\begin{array}{l}
    \varphi_\varepsilon''+c A_\varepsilon'+\varepsilon^2(2A_\varepsilon' \varphi_\varepsilon'+A_\varepsilon\varphi_\varepsilon'')=0, \\
    -c\varphi_\varepsilon'+2f'(1) A_\varepsilon+\varepsilon^2\Big( A_\varepsilon^2\big(3f'(1)+2f''(1)\big)+A_\varepsilon''-cA_\varepsilon\varphi_\varepsilon'-(\varphi_\varepsilon')^2\Big)=O(\varepsilon^4),
    \end{array}
\right.
\end{equation}

where $h(\varepsilon)=O(\varepsilon^\alpha)$ means that we can find $C>0$ that does not depend on $\varepsilon$ such that $|h(\varepsilon)|\leq C\varepsilon^\alpha$.

For the speed $c=\sqrt{c_s^2 - \varepsilon^2}$ taken in the limit $c\rightarrow c_s$ ($\varepsilon\rightarrow 0$), \eqref{deuxieme systeme TWc en remplacant u par l'ansatz u_eps} leads to

\begin{equation}\label{-c_s phi' eps + 2f'(1)A eps = O(eps^2)}
    \left\{
\begin{array}{l}
    \varphi_\varepsilon''+c_s A_\varepsilon'=O(\varepsilon^2), \\
    \varphi_\varepsilon'+c_s A_\varepsilon=O(\varepsilon^2).
    \end{array}
\right.
\end{equation}

Differentiating the second equality of \eqref{deuxieme systeme TWc en remplacant u par l'ansatz u_eps} and adding the first one, multiplied by $c$, we obtain 

\begin{equation*}
    -A_\varepsilon'+\big(3f'(1)+2f''(1)\big)2A_\varepsilon A_\varepsilon ' + A_\varepsilon'''+c A_\varepsilon ' \varphi_\varepsilon'-2\varphi_\varepsilon ' \varphi_\varepsilon''=O(\varepsilon^2).
\end{equation*}

Using \eqref{-c_s phi' eps + 2f'(1)A eps = O(eps^2)} and that $c=c_s + O(\varepsilon^2)$, we get 
\begin{equation*}
    -A_\varepsilon'+\big(6f'(1)+2f''(1)\big)2A_\varepsilon A_\varepsilon'+A_\varepsilon'''=0.
\end{equation*}

Still in the limit $c\rightarrow c_s$, assuming that $A_\varepsilon$ and $\varphi_\varepsilon$ converge to some functions $A$ and $\varphi$, we obtain

\begin{equation*}
    -A'+\big(6f'(1) + 2 f''(1)\big)2A A'+A'''=0,
\end{equation*}

and assuming that $A,A',A'' \underset{|x|\rightarrow +\ii}{\longrightarrow}0$, we are led to

\begin{equation}\label{kdv equation avec k=2f''(1)+6f'(1)}
    -A+A''+kA^2=0\quad\text{with }k:=2f''(1)+6f'(1).
\end{equation}

Observe that this is the integrated version of the $(KdV)$ equation that appears in~\cite{Chiron9} with $$\Gamma=\dfrac{-2k}{c_s^2}.$$

The condition \eqref{condition suffisante pour la stabilité orbitale sur f''(1)+6f'(1)>0} is equivalent to the fact that $k\neq 0$. Our choice for $A_\varepsilon$ corresponds to a soliton for the (KdV) equation \eqref{kdv equation avec k=2f''(1)+6f'(1)}. Its expression is given explicitly by

\begin{equation}\label{expression de A, solution de KDV avec k}
    A(x):=\dfrac{3}{2k}\sech^2(\dfrac{x}{2}).
\end{equation}

Going back to \eqref{-c_s phi' eps + 2f'(1)A eps = O(eps^2)} in the limit $\varepsilon \rightarrow 0$ gives $\varphi'=-c_s A$, so that we are led to the choice

\begin{equation}\label{expression de phi (kdv), avec k}
    \varphi(x):=\dfrac{-3c_s}{k}\tanh(\dfrac{x}{2}).
\end{equation}

\begin{lem}
Suppose that $k\neq 0$. Let $v_\varepsilon=\big(1+\varepsilon^2 A(\varepsilon .)\big) e^{i\varepsilon\varphi(\varepsilon .)}$, where $A$ and $\varphi$ are given by \eqref{expression de A, solution de KDV avec k} and \eqref{expression de phi (kdv), avec k}. Then

\begin{equation}\label{lem: estimée energie et moment en eps^3 et eps^5 +O(eps^7)}
    E(v_\varepsilon)=\dfrac{6c_s^2}{k^2}\varepsilon^3-\dfrac{18}{5k^3}\big(f''(1)+5f'(1)\big)\varepsilon^5 + O(\varepsilon^7)\quad\text{and}\quad p(v_\varepsilon)=\dfrac{6 c_s}{k^2}\Big( \varepsilon^3+\dfrac{3}{5k}\varepsilon^5\Big).
\end{equation}
\end{lem}

\begin{rem}
These asymptotic expansions were already computed in Theorem 4 in~\cite{Chiron7} using the integral formulas of the energy and of the momentum.
\end{rem}

\begin{proof}
Bearing in mind that $\varphi'=-c_s A$ and using the following identities
\begin{equation*}
    \int_\R \sech^4(\dfrac{x}{2})dx=\dfrac{8}{3},\quad\int_\R \sech^6(\dfrac{x}{2})dx=\dfrac{32}{15}\quad\text{and}\quad\int_\R \sech^4(\dfrac{x}{2})\tanh^2(\dfrac{x}{2})dx = \dfrac{8}{15},
\end{equation*}

we first compute the momentum

\begin{align*}
    p(v_\varepsilon) &= -\dfrac{1}{2}\int_\R \big(2\varepsilon^2 A(\varepsilon.)+\varepsilon^4 A(\varepsilon .)^2\big)\varepsilon^2 \varphi'(\varepsilon .)\\
    &=\dfrac{6 c_s}{k^2}\varepsilon^3 +\dfrac{18 c_s}{5k^3}\varepsilon^5.
\end{align*}

Now, we compute the kinetic energy

\begin{align*}
    E_k(v_\varepsilon) &=\dfrac{1}{2}\int_\R \big( 1+\varepsilon^2 A(\varepsilon .)\big)^2\varepsilon^4\varphi'(\varepsilon .)^2+\varepsilon^6 A'(\varepsilon .)^2\\
    &=\dfrac{3c_s^2}{k^2}\varepsilon^3+\dfrac{36c_s^2}{5k^3}\varepsilon^5+\dfrac{3}{5k^2}\varepsilon^5+O(\varepsilon^7).
\end{align*}

For the potential energy, we first write the Taylor expansion of $F$. We claim that if $z\in L^4(\R)\cap L^8(\R)$, then there exists $C >0$, only depending on $f$, such that

\begin{equation}\label{lem: developpement de taylor de F avec R_eps(z)}
F\big((1+\varepsilon^2 z)^2\big)=c_s^2 z^2 \varepsilon^4 + \Big(c_s^2-\dfrac{4 f''(1)}{3}\Big)z^3 \varepsilon^6 + \varepsilon^8 \widetilde{R}_\varepsilon(z),
\end{equation}

with
\begin{equation*}
    \normLun{\widetilde{R}_\varepsilon(z)} \leq C \sum_{k=4}^8 \Vert z\Vert_{L^k}^k.
\end{equation*}

Indeed, we compute the fourth order Taylor expansion of the function $x\mapsto F(1+x)$. For $x >0$, 

\begin{equation*}
    F(1+x)=-\dfrac{f'(1)}{2}x^2-\dfrac{f''(1)}{6}x^3-\dfrac{x^4}{6}\int_0^1 (1-t)^3 f'''(1+tx)dt.
\end{equation*}

Let $z$ be a function in $L^4(\R)\cap L^8(\R)$. Replacing $x$ by $2\varepsilon^2 z + \varepsilon^4 z^2$, and invoking the previous expansion,

\begin{align*}
    F\big((1+\varepsilon^2 z)^2\big)&=c_s^2 z^2 \varepsilon^4+\Big(c_s^2-\dfrac{4 f''(1)}{3}\Big)z^3 \varepsilon^6 + \varepsilon^8 \widetilde{R}_\varepsilon(z)
\end{align*}

with

\begin{align*}
     \widetilde{R}_\varepsilon(z)=-\dfrac{f''(1)}{6}(z^4 + 6\varepsilon^2z^5 + \varepsilon^4 z^6)-\dfrac{(2z+\varepsilon^2 z^2)^4}{6}\int_0^1 (1-t)^3 f'''(1+t\big(2\varepsilon^2 z + \varepsilon^4 z^2)\big))dt.
\end{align*}

We set $z=A(\varepsilon .)\in H^1(\R)\hookrightarrow L^p(\R)$ for any $2\leq p\leq +\ii$. In particular, this is true with $p=+\ii$, then $1+t(2\varepsilon^2 z+\varepsilon^4 z^2)\in L^{\ii}(\R)$. By hypothesis, $f'''$ is continuous, that is why there exists $M>0$ independent of $x$ and $\varepsilon$ such that $$\Big|\int_0^1 (1-t)^3 f'''(1+t\big(2\varepsilon^2 z + \varepsilon^4 z^2)\big))dt\Big|\leq M.$$

Now using~\eqref{lem: developpement de taylor de F avec R_eps(z)} and the Sobolev embeddings described just above, we are led to

\begin{align*}
    E_p(v_\varepsilon) &=\dfrac{1}{2} \int_\R \bigg( c_s^2 A(\varepsilon .)^2 \varepsilon^4 + \Big(c_s^2-\dfrac{4 f''(1)}{3}\Big)A(\varepsilon .)^3 \varepsilon^6 + \varepsilon^8 \widetilde{R}_\varepsilon\big(A(\varepsilon .)\big)\bigg)\\
    &= \dfrac{3c_s^2}{k^2}\varepsilon^3 + \Big(c_s^2-\dfrac{4 f''(1)}{3}\Big)\dfrac{18}{5k^3}\varepsilon^5 + O(\varepsilon^7).
\end{align*}

Adding both $E_p(v_\varepsilon)$ and $E_k(v_\varepsilon)$ expressions, we obtain

\begin{align*}
    E(v_\varepsilon)&=\dfrac{6c_s^2}{k^2}\varepsilon^3+\dfrac{3k+54c_s^2-24f''(1)}{5k^3}\varepsilon^5 + O(\varepsilon^7)\\
    &=\dfrac{6c_s^2}{k^2}\varepsilon^3-\dfrac{18}{5k^3}\big(f''(1)+5f'(1)\big)\varepsilon^5 + O(\varepsilon^7).
\end{align*}
\end{proof}

\begin{proof}[Proof of Proposition \ref{prop inegalité a droite de Emin par des puissances de q}.]
For $\mathfrak{q}$ small enough, we can parametrize $\mathfrak{q}$ as a function of $\varepsilon$ as 

\begin{equation*}
    \mathfrak{q}_\varepsilon = \dfrac{6 c_s}{k^2}\Big( \varepsilon^3+\dfrac{3}{5k}\varepsilon^5\Big).
\end{equation*}

Indeed, no matter how we take $k$, as long as it is non-zero, $\varepsilon\mapsto\dfrac{6 c_s}{k^2}\Big( \varepsilon^3+\dfrac{3}{5k}\varepsilon^5\Big)$ is a strictly increasing function of $\varepsilon\in (0, \sqrt{|k|})$. For simplicity, we set 

\begin{equation*}
    \mathfrak{s}_\varepsilon :=\dfrac{k^2}{6c_s}\mathfrak{q}_\varepsilon=\varepsilon^3 + \dfrac{3}{5k}\varepsilon^5.
\end{equation*}

We address the case $k <0$ for which we assume that $\varepsilon$ is in the interval $(0, \sqrt{|k|})$. We have

\begin{equation}\label{prop: expression de s_eps dans le cas k<0}
    \mathfrak{s}_\varepsilon =\varepsilon^3 - \dfrac{3}{5|k|}\varepsilon^5,
\end{equation}

so that 

\begin{equation}\label{prop: inégalité à gauche et à droite de s_eps en eps^3}
    \dfrac{2}{5}\varepsilon^3\leq \mathfrak{s}_\varepsilon \leq \varepsilon^3 .
\end{equation}

Applying the Taylor-Lagrange theorem, and observing that $\mathfrak{s}_\varepsilon \geq \dfrac{2}{5|k|}\varepsilon^5$ we find $\mathfrak{p}_{\varepsilon,k}\in (\mathfrak{s}_\varepsilon,7\mathfrak{s}_\varepsilon)$, such that

\begin{equation*}
    \varepsilon^5 = \gs_\varepsilon^{\frac{5}{3}}+\dfrac{\varepsilon^5}{|k|}\mathfrak{p}_{\varepsilon,k}^{\frac{2}{3}}.
\end{equation*}

Using again (\ref{prop: inégalité à gauche et à droite de s_eps en eps^3}), we conclude that 

\begin{equation*}
    \varepsilon^5 =\gs_\varepsilon^{\frac{5}{3}}+ O(\gs_\varepsilon^{\frac{7}{3}})=\Big(\dfrac{k^2}{6c_s}\Big)^{\frac{5}{3}}\gq_\varepsilon^{\frac{5}{3}}+ O(\gq_\varepsilon^{\frac{7}{3}}).
\end{equation*}

Combining this asymptotics with the expression of $E(v_\varepsilon)$ in \eqref{lem: estimée energie et moment en eps^3 et eps^5 +O(eps^7)}, \eqref{prop: expression de s_eps dans le cas k<0} and \eqref{prop: inégalité à gauche et à droite de s_eps en eps^3}, we get

\begin{align*}
    E(v_\varepsilon)&=\dfrac{6c_s^2}{k^2}\Big(\gs_\varepsilon + \dfrac{3}{5|k|}\varepsilon^5\Big) - \dfrac{18}{5k^3}\big( f''(1)+5 f'(1)\big)\varepsilon^5 + O(\gq_\varepsilon^{\frac{7}{3}})\\
    &=c_s \mathfrak{q}_\varepsilon -\dfrac{9 (k^2)^{\frac{2}{3}}}{5(6c_s)^{\frac{5}{3}}}\gq_\varepsilon^{\frac{5}{3}}+ O(\gq_\varepsilon^{\frac{7}{3}})
\end{align*}

\end{proof}

Now let us handle the case $k>0$. By the same argument, rather noticing that $\varepsilon^3\leq \gs_\varepsilon \leq \frac{8}{5}\varepsilon^3$ and that $\gs_\varepsilon \geq \frac{4}{5k}\varepsilon^5$, we deduce the existence of a similar number $\gp_{\varepsilon,k}\in (\frac{1}{4}\gs_\varepsilon,\gs_\varepsilon)$ and this leads to the same asymptotic expression of $E(v_\varepsilon)$. In both cases ($k<0$ and $k>0$), we set $K_1=\frac{9 (k^2)^{\frac{2}{3}}}{5(6c_s)^{\frac{5}{3}}}\gq_\varepsilon^{\frac{5}{3}}$ and the number $K_2>0$ is provided by $O(\gq_\varepsilon)$. Hence Proposition \ref{prop inegalité a droite de Emin par des puissances de q} is proved.

\begin{rem}
When $k= 0$, we lose the nonlinear effects and the (KdV) limit equation turns into a linear equation. According to the article by D. Chiron~\cite{Chiron9}, it is natural to change the scaling and use a different ansatz. We obtain a modified (KdV) equation where an analogue constant $\Gamma'$ involves the quantities $c_s$ and $f'''(1)$. There exist solutions whenever $\Gamma' <0$, and then repeating the same method as above, we may find the same type of estimates for $\Emin$. However, we can expect different exponents than those exhibited in~\eqref{prop: inegalité a droite de Emin par des puissances de q}, in view of the difference between~\eqref{lem: estimée energie et moment en eps^3 et eps^5 +O(eps^7)} and the analogue estimates in Theorem~5 in~\cite{Chiron7}.
\end{rem}

\begin{cor}\label{corollaire Emin is striclty sub additive}
The function $E_{\min}$ is strictly subadditive on $\R_+$ and for $\gq >0$, it satisfies $E_{\min}(\gq) < c_s \gq$.
\end{cor}

\begin{proof}
We use a general result on continuous concave functions that vanish in $0$, that is the following lemma.
\begin{lem}[\cite{BetGrSa1},\cite{ChirMar2}] \label{appendix lemme pour fonction concave, soit sous additive soit lineaire}
Let $f:[0,+\ii[\rightarrow \R$ be a continuous concave fonction, with $f(0)=0$ and owning a finite right derivative at the origin

\begin{equation*}
    a:=\lim_{x\rightarrow 0^+}\dfrac{f(x)}{x}.    
\end{equation*}

Then for any $\mathfrak{s}>0$, the following alternative holds:
\begin{itemize}
    \item $f$ is linear on $[0,\mathfrak{s}]$, with slope $a$, or
    \item $f$ is strictly subadditive.
\end{itemize}
\end{lem}

Combining this lemma with estimate  \eqref{prop: inegalité a droite de Emin par des puissances de q} implies that $\Emin$ cannot be linear near zero and then it is strictly subadditive on $\R_+$. Moreover, if $\gq >0$,

\begin{align*}
    \Emin (\gq)  = \Emin \Big(\dfrac{\gq}{2}+\dfrac{\gq}{2}\Big) &< 2 \Emin\Big(\dfrac{\gq}{2}\Big)\quad\text{by strict subadditivity,}\\
    & \leq c_s \gq\quad\text{by estimate \eqref{lem: Eminq leq cs q}.}
\end{align*}
\end{proof}

\section{Characteristics of the kink solution}\label{section charcteristics of the kink solution}
\begin{prop}\label{prop: la solution de tw0 est unique si int_0^1 f neq 0}
Let $f$ be a continuous function such that \eqref{1ere hypothèse de croissance sur F minorant intermediaire} holds. Then there exists a solution $\mathfrak{v}_0$ of \eqref{TWC} with $c=0$ and such that $|\gv_0(x)|\underset{|x|\rightarrow +\ii}{\longrightarrow}1$. This solution is unique up to a constant phase shift and a translation.
\end{prop}

\begin{proof}
The proof is based on Theorem~\ref{existence et unicité des travelling wave si condition du zéro}. We set $\xi_0=1$ and verify that $\mathcal{N}_0(\xi_0)=0$. In addition, using assumption \eqref{1ere hypothèse de croissance sur F minorant intermediaire}, we get $\mathcal{N}_0(\xi)=-4(1-\xi)F(1-\xi)\leq -4\lambda(1-\xi)\xi^2<0$ for $\xi\in (0,\xi_0)$.

Using again assumption \eqref{1ere hypothèse de croissance sur F minorant intermediaire}, we have \begin{equation*}
    \mathcal{N}_0'(\xi_0)=4F(0) \geq 4\lambda >0,
\end{equation*}

so that the assumptions of Theorem~\ref{existence et unicité des travelling wave si condition du zéro} are satisfied. The existence and uniqueness are consequences of Theorem~\ref{existence et unicité des travelling wave si condition du zéro}.
\end{proof}

Let us emphasize the specific role of the black soliton $\mathfrak{v}_0$.

\begin{lem}\label{lemme l'inf de l'énergie dont la fonction s'annule est atteint en u_0}
We have \begin{equation*}
    E(\mathfrak{v}_0)=\inf\big\lbrace E(v)\big| v\in H^1_{\loc}(\R), \inf_\R |v|=0\big\rbrace .
\end{equation*}

In particular, if $E(v) < E(\mathfrak{v}_0)$, then \begin{equation*}
    \inf_\R |v| >0.
\end{equation*}
\end{lem}

\begin{proof}
We consider a minimizing sequence $(u_n)$ for the problem

\begin{equation*}
    \mathcal{E}_0 :=\inf\Big\lbrace \int_0^{+\ii}e(v) \big| v\in H^1_{\loc}([0,+\ii[), v(0)=0\Big\rbrace.
\end{equation*}

~~\\

The sequence $(u_n')$ is bounded in $L^2([0,+\ii))$ with respect to $n$. We also have $u_n(0)=0$ so that we obtain

\begin{equation*}
    \int_{0}^R |u_n(x)|^2 dx\leq \dfrac{R^2}{2} \Vert u_n'\Vert_{L^2([0,+\ii[)}^2,
\end{equation*}

which shows that $(u_n)$ is bounded in $H^1_{\loc}([0,+\ii))$. Hence, using the Rellich compactness theorem, we can assume, up to a subsequence, that $(u_n)$ tends strongly in $L_{\loc}^\ii ([0,+\ii))$ to a function $u$. Still up to a subsequence, we can assume that $(u_n')$ tends weakly to $u'$ in $L^2([0,+\ii))$.\\

By Fatou's lemma and by the weak convergence of $(u'_n)$, we are led to

\begin{equation*}
    \int_0^{+\ii} e(u)=\dfrac{1}{2}\int_0^{+\ii}(u')^2+\dfrac{1}{2}\int_0^{+\ii}\liminf_{n\rightarrow +\ii}F(|u_n|^2)\leq \liminf_{n\rightarrow +\ii}\int_0^{+\ii} e(u_n).
\end{equation*}

Thus the infimum is achieved by the function $u$ which is then critical for the Ginzburg-Landau energy i.e. 

\begin{equation*}
    0=\nabla E(u)=- u'' - u f(|u|^2).
\end{equation*}

We can extend $u$ to an odd function. Indeed, since $(u_n)$ converges in $L^{\ii}_{\loc}([0,+\ii))$, the limit function $u$ is continuous on $\R_+$, and then the convergence is pointwise. Then the fact that $u_n(0)=0$ implies $u(0)=0$. Moreover, by the equation, $u$ actually lies in $\mathcal{C}^2\big([0,+\ii)\big)$, so setting $v:=-u(-.)$, it satisfies $v(0)=u(0)$ and $v'(0)=u'(0)$ so that $v$ solves the same Cauchy problem as $u$ in $0$, hence the equality $u=v$. Moreover, it proves that $u'(0)\neq 0$, otherwise we would have $u\equiv 0$ which contradicts the non-vanishing property at infinity. Proposition~\ref{prop: la solution de tw0 est unique si int_0^1 f neq 0} then implies that $u=\mathfrak{v}_0$ by uniqueness (up to a phase shift).

For each function $v\in H^1([0,+\ii[)$ with finite energy and that vanishes (we can assume that it vanishes in $0$), we obtain by minimality of $\mathfrak{v}_0$ 

\begin{equation*}
    \int_0^{+\ii} e(v)\geq \int_0^{+\ii} e(\mathfrak{v}_0).
\end{equation*}

The same argument holds on $(-\ii,0]$ by oddness and then $E(v)\geq E(\mathfrak{v}_0)$. Thus \begin{equation*}
    E(\mathfrak{v}_0)=\inf\big\lbrace E(v)\big| v\in H^1_{\loc}(\R),\inf_\R |v|=0\big\rbrace.
\end{equation*}

In particular, if $E(v)<E(\mathfrak{v}_0)$, then $\inf_\R |v|>0$, which shows that $v$ never vanishes.
\end{proof}

\begin{rem}
We have proven in addition that the kink solution is odd and vanishes at $x=0$.
\end{rem}

A consequence of this property is the following.
\begin{prop}\label{prop Emin(q) leq Emin(pi/2) = E(v_0)}
For $\gq\geq\gq_*$, we have \begin{equation*}\label{prop: égalité Emin(pi/2) = Ecyrilique min(pi/2)}
    E_{\min}(\gq)=E(\mathfrak{v}_0).
\end{equation*}
\end{prop}

\begin{proof}
We prove that, for any $\gq\geq \gq_*$, there exists a sequence $(w_n)\in \Nenergyset^\N$ such that 

\begin{equation*}
    E(w_n)\ntend E(\mathfrak{v}_0)\quad\text{and }\quad p(w_n)=\gq .
\end{equation*}

We set 

\begin{equation*}
    w_n(x) = \left\{
\begin{array}{l}
    \Big|\gv_0\Big(\dfrac{1}{n}\Big)\Big|e^{i\psi_n(x)}\quad\text{if }|x|\leq \dfrac{1}{n}, \\
    |\gv_0(x)|e^{i\psi_n(x)}\quad\text{if }|x|\geq \dfrac{1}{n},
\end{array}
\right.
\end{equation*}

with 

$$\psi_n(x) :=\left\{
\begin{array}{l}
    0\quad\text{if }x\geq\frac{1}{n}, \\
    q_n(nx-1)\quad\text{if }|x|\leq  \dfrac{1}{n},\\
    -2q_n\quad\text{if }x\leq -\frac{1}{n},
\end{array}
\right.$$ 

and

\begin{equation*}
    q_n= \dfrac{\gq}{1-\Big|\gv_0\Big(\dfrac{1}{n}\Big)\Big|^2}.
\end{equation*}

We verify that $w_n\in\Nenergyset$ for all $n\in\N^*$, with

\begin{equation*}
    p(w_n)=\dfrac{1}{2}\int_{-\frac{1}{n}}^{\frac{1}{n}}\bigg( 1- \Big|\gv_0\Big(\dfrac{1}{n}\Big)\Big|^2\bigg)\psi_n' = \gq .
\end{equation*}

Moreover, we compute
\begin{align*}
    E(w_n) & = \dfrac{1}{2}\int_\R |w_n'|^2 + \dfrac{1}{2}\int_\R F(|w_n|^2)\\ 
    & = \dfrac{1}{2}\int_{|x|\leq\frac{1}{n}} (\psi_n')^2 \Big|\gv_0\Big(\dfrac{1}{n}\Big)\Big|^2 +\dfrac{1}{2}\int_{|x|\geq\frac{1}{n}} (|\gv_0|')^2 +\dfrac{1}{2}\bigg( \int_{|x|\leq\frac{1}{n}} F(|w_n|^2)+\int_{|x|\geq \frac{1}{n}}F(|w_n|^2)\bigg)\\
    & = n q_n^2 \Big|\gv_0\Big(\dfrac{1}{n}\Big)\Big|^2 + \dfrac{1}{2}\int_{|x|\geq \frac{1}{n}}|\gv_0'|^2 +  \dfrac{1}{2}\int_{|x|\leq\frac{1}{n}}F\Big(\Big|\gv_0\Big(\dfrac{1}{n}\Big)\Big|^2\Big)+\dfrac{1}{2}\int_{|x|\geq\frac{1}{n}}F(|\gv_0|^2).
\end{align*}

By continuity of $F$, we obtain that 

\begin{equation*}
    \int_{-\frac{1}{n}}^{\frac{1}{n}}F\Big(\Big|\gv_0\Big(\dfrac{1}{n}\Big)\Big|^2\Big)\ntend 0.
\end{equation*}

On the other hand, since $q_n\longrightarrow \gq$ as $n\rightarrow +\ii$, and by differentiability of $\gv_0$ at $x=0$,

\begin{equation*}
    n q_n^2 \Big|\gv_0\Big(\dfrac{1}{n}\Big)\Big|^2 \leq \dfrac{C^2 q_n^2}{n}\ntend 0.
\end{equation*}

Finally, we get the limit

$$E(w_n)\ntend E(\gv_0)\quad\text{with }p(w_n)=\gq.$$

This shows that

$$\Emin(\gq)\leq E(\gv_0).$$

By contradiction, let us assume that $\Emin(\gq_*)<E(\gv_0)$, then by continuity, there exists $\varepsilon >0$ such that $\Emin(\gq_*+\varepsilon)<E(\gv_0)$. In view of Lemma~\ref{lemme l'inf de l'énergie dont la fonction s'annule est atteint en u_0}, this contradicts the definition of $\gq_*$. Finally, since $\Emin$ is nondecreasing on $\R_+$, we conclude that $\Emin(\gq)=\Emin(\gq_*)$ for all $\gq\geq \gq_*$.
\end{proof}

We can now prove that $\Emin$ is stricty increasing on $[0,\gq_*]$.
\begin{proof}[Proof of Proposition~\ref{lemme Emin est strictement croissante}]
We already know by Lemma~\ref{lemme Emin est croissante} that $\Emin$ is nondecreasing. Now, if it is not strictly increasing on $[0,\gq_*)$, then there exists $0\leq a < b <\gq_*$ such that $\Emin$ is constant on $[a,b]$. By Proposition~\ref{proposition Emin est concave}, $\Emin$ is also concave on $\R_+$ which implies that $\Emin$ is constant on $[a,+\ii)$ and then that $\Emin (a)=\Emin(\gq_*)$. Then, by Proposition~\ref{prop: égalité Emin(pi/2) = Ecyrilique min(pi/2)}, we have $\Emin(a)=E(\gv_0)$, therefore $\gq_*\leq a$ which brings a contradiction.
\end{proof}

We also try to give a rigorous meaning to $p(\mathfrak{v}_0)$ which is not so immediate because of the vanishing property of the black soliton. Let us first recall how to define properly the momentum of a function in $\energyset$.

\begin{lem}[\cite{BeGrSaS1}]\label{lemme: formule du moment [p] sur chi^1 tilde}
Let $v\in \energyset$. Then the limit 

\begin{equation*}
    [p](v)=\lim_{R\rightarrow +\ii}\Big( \int_{-R}^R \langle iv',v\rangle +\dfrac{1}{2}\big( \arg v(R) - \arg v(-R)\big)\Big) \mod \pi
\end{equation*}
exists. Moreover, if $v\in\mathcal{N}\energyset$, then 
\begin{equation*}
    [p](v)=p(v) \mod\pi .
\end{equation*}

\end{lem}

\begin{prop}\label{proposition untwisted momentum de v_0}
    We have 
    \begin{equation*}
        [p](\gv_0)=\dfrac{\pi}{2}\mod\pi.
    \end{equation*}
\end{prop}

For $c\in\R$, a solution $\gv_c$ to~\eqref{TWC} satisfies the general formula

\begin{equation*}
    \langle i\mathfrak{v}_c,\mathfrak{v}_c'\rangle = \dfrac{c}{2}\eta_c .
\end{equation*}

Indeed, writing $v_1$ (respectively $v_2$) the real part (respectively the imaginary part) of $\gv_c$, we obtain simultaneously 

$$\langle i\mathfrak{v}_c,\mathfrak{v}_c'\rangle = v_1v_2' - v_2 v_1',$$

and 

$$\left\{
\begin{array}{l}
    v_1 '' - cv_2' +v_1f(|v|^2)=0, \\
    v_2'' + cv_1 + v_2f(|v|^2)=0. \\
\end{array}
\right.$$

Hence, by multiplying the first line by $v_2$ and the second one by $v_1$, we get $(v_1 v_2' - v_2 v_1')'=\frac{c}{2}\eta_c'$ then $\langle i\mathfrak{v}_c,\mathfrak{v}_c'\rangle = \frac{c}{2}\eta_c$.
Therefore $\langle i\mathfrak{v}_0,\mathfrak{v}_0'\rangle =0$\ \footnote{This can also be seen by the fact that $\gv_0$ is real-valued}. In addition, due to the oddness property, the limit at infinity and the fact that $\gv_0$ is real-valued, we obtain either that $\gv_0(R)$ is close to $1$ and $\gv_0(-R)$ to $-1$, or the contrary, for $R$ large enough. In other words, either $\arg\big(\gv_0(R)\big)=0$ and $\arg\big(\gv_0(-R)\big)=\pi$, or $\arg\big(\gv_0(R)\big)=\pi$ and $\arg\big(\gv_0(-R)\big)=0$. In both cases, $\arg \mathfrak{v}_0(R)-\arg \mathfrak{v}_0(-R)=-\pi$ (respectively $\pi$), so that

$$[p](\mathfrak{v}_0) =\frac{\pi}{2}\mod\pi.$$

\section{Minimization of the energy at fixed momentum}\label{section méthode variationnelle, minimisation de l'énergie à moment fixé}
We now handle the construction of minimizing travelling waves for a general nonlinearity with non-vanishing condition at infinity. Whenever $\gp\in (0,\gq_*)$ avoids the critical momentum of the black soliton, we prove that there exist minimizers for $\Emin(\gp)$. We need to analyze minimizing sequences, and for dealing with orbital stability, we also consider pseudo-minimizing sequences. These are sequences $(u_n)\in\Nenergyset^\N$ satisfying

\begin{equation}\label{etre une suite pseudo minimisante}
    p(u_n)\ntend\mathfrak{p}\quad\text{and}\quad E(u_n)\ntend \Emin (\mathfrak{p}),
\end{equation}

Here, we fix $\gp$ such that

\begin{equation}\tag{$H_{\gq_*}$}
    \left\{
\begin{array}{l}
    \quad\gp\in (0,\gq_*)\\
    \quad\quad\text{ or }\quad\\
    \gp=\gq_*\notin \frac{\pi}{2}+\pi\Z. \\ 
\end{array}
\right.
\end{equation}

In order to prove Theorem~\ref{theoreme de continuité de l'energie dans le cas d'une non linéarité générale} , we settle a concentration-compactness argument. This argument relies on separate lemmas that we now present. Henceforth, we assume that \eqref{hypothèse de croissance sur F minorant intermediaire}, \eqref{hypothèse de croissance sur F majorant} and \eqref{condition suffisante pour la stabilité orbitale sur f''(1)+6f'(1)>0} hold.

\begin{lem}\label{lemme controle de v pres et loin de 1}
Let $E >0$ and $\delta \in (0,1)$ be given. Then there exists $l_0\in\N^*$, depending only on the previous quantities, such that the following property holds. Given any map $v\in H^1_\loc(\R)$ satisfying $E(v)\leq E$, either 

\begin{equation*}
    \big| 1-|v(x)|\big| <\delta_0\quad\forall x \in\R,
\end{equation*}

or there exists $l\leq l_0$ points $x_1,...,x_l$ such that 

\begin{equation*}
    \big| 1-|v(x_i)|\big| \geq\delta_0\quad\forall i\in\{ 1,...,l\},
\end{equation*}

and 
 \begin{equation*}
    \big| 1-|v(x)|\big| <\delta_0\quad\forall x \in \R\setminus\bigcup_{i=1}^l [x_i-1,x_i + 1 ].
\end{equation*}

\end{lem}

\begin{proof}
We set $$\mathcal{A}:= \big\lbrace z\in\R\big|\ \big| 1-|v(z)|\big|\geq \delta_0\big\rbrace,$$ 

and assume that this set is not empty. Consider also the covering with the intervals $I_i:=[i-\frac{1}{2},i+\frac{1}{2}]$. We claim that if $\mathcal{A}\cap I_i \neq \varnothing$, then \begin{equation}\label{int_widetilde I_i e(v) geq mu_0}
    \int_{\widetilde{I}_i}e(v)\geq \mu_0,
\end{equation}

where $\widetilde{I}_i=[i-1,i+1]$ and $\mu_0$ is some positive constant. To prove the claim, we first notice that 
for any $(x,y)\in\R^2$, we have

\begin{equation*}\label{|v(x)-v(y)| leq ...}
    \big| v(x)-v(y)\big| \leq \Vert v'\Vert \sqrt{|x-y|}\leq \sqrt{2E|x-y|}.
\end{equation*}

Therefore, there exists $r>0$ such that, if $z\in\mathcal{A}$, then for all $y\in [z-r,z+r]$, 

$$\big| 1-|v(y)|\big| \geq \dfrac{\delta_0}{2}.$$

Choosing $r_0:=\min(r,\frac{1}{2})$ and invoking \eqref{1ere hypothèse de croissance sur F minorant intermediaire}, we are led to \begin{align*}
    \int_{z-r_0}^{z+r_0} e(v) & \geq \dfrac{1}{2}\int_{z-r_0}^{z+r_0}F(|v|^2)\\
    & \geq \dfrac{\lambda}{2}\int_{z-r_0}^{z+r_0}(1-|v|^2)^2\\
    & \geq \dfrac{\lambda}{2}\int_{z-r_0}^{z+r_0}(1-|v|)^2 \geq \mu_0 :=\dfrac{\lambda r_0 \delta_0^2}{4}.
\end{align*}

In particular, if $z\in I_i\cap \mathcal{A}$ for some $i\in\N$, then $[z-r_0,z+r_0]\subset \widetilde{I}_i$, and claim \eqref{int_widetilde I_i e(v) geq mu_0} follows. To conclude the proof, we notice that 

\begin{equation*}
    \sum_{i\in\N}\int_{\widetilde{I}_i}e(v)=2E(v)\leq 2E ,
\end{equation*}

so that, in view of \eqref{int_widetilde I_i e(v) geq mu_0}, 

\begin{equation*}
    l\mu_0\leq 2E,
\end{equation*}

where $l:=\#\lbrace i\in\N | I_i\cap\mathcal{A}\neq\varnothing\rbrace$. The conclusion follows choosing $l_0$ as the greatest integer below $\frac{2E}{\mu_0}$ and choosing $x_i\in I_i\cap\mathcal{A}$ for any $i\in\N$ such that $I_i\cap\mathcal{A}\neq \varnothing$.
\end{proof}

We also need the following construction.\\

\begin{lem}\label{lemme equivalent du lemme 6 dans existance and travelling Begrasa 2}
Let $|\gq| < \frac{1}{32}$ and $0\leq \mu\leq \frac{1}{4}$. There exist some number $\ell >1$, a map $w=|w|e^{i\varphi}\in H^1([0,\ell])$, and a number $C>0$ depending only on $f$, such that 

\begin{equation*}
    w(0)=w(\ell)\quad\text{and}\quad \big| 1- |w(0)|\big|\leq\mu,
\end{equation*}

\begin{equation*}
    \gq=\dfrac{1}{2}\int_0^\ell (1-|w|^2)\varphi',
\end{equation*}

and \begin{equation}\label{lemme  E(w) leq C |q|}
    E(w)\leq C |\gq|.
\end{equation}
\end{lem}

\begin{proof}
As in Lemma 6 in~\cite{BetGrSa2}, we construct for $\lambda >0$, functions $f_\lambda=\frac{1}{\lambda}f(\frac{.}{\lambda})$ and $\varphi_\lambda=\varphi(\frac{.}{\lambda})$ with

\begin{equation*}
    f(s)=\left\{
\begin{array}{l}
    s \quad\text{if }s\in [0,\frac{1}{2}],\\
    1-s\quad \text{if }s\in [\frac{1}{2},1 ],\\
    0\quad \text{if }s\in [1,2],
    \end{array}
\right. \quad\text{and}\quad  \varphi(s)=\left\{
\begin{array}{l}
    s \quad\text{if }s\in [0,1 ],\\
    2-s\quad \text{if }s\in [1,2].\\
    \end{array}
\right.
\end{equation*}

They satisfy
\begin{equation*}
    |f_\lambda|\leq \dfrac{1}{2\lambda},|\varphi'_\lambda|=\dfrac{1}{\lambda}, f_\lambda(0)=f_\lambda(2\lambda)=0,\varphi_\lambda(0)=\varphi_\lambda(2\lambda)=0.
\end{equation*}

We choose $\lambda =\frac{1}{8|\gq|}$, so that $\frac{1}{\lambda}\leq\frac{1}{4}$ and then $f_\lambda\leq \frac{1}{8}$. Introduce a parameter $\delta \in [0,\frac{1}{2}]$ such that $1-\sqrt{1-\delta}\leq \mu$, and consider the function 

\begin{equation*}
    \rho_{\lambda,\delta}=\sqrt{1-\delta-f_\lambda},
\end{equation*}

so that $\rho_{\lambda,\delta}^2$ is bounded in the interval $[0,2\lambda]$. It follows from the special choice of parameters $\lambda$ and $\delta$ that 

\begin{equation*}
    |\gq|=\dfrac{1}{2}\int_0^{2\lambda} (1-\rho_{\lambda,\delta}^2)\varphi_\lambda' \quad\text{and}\quad |\delta + f_\lambda|\leq 1.
\end{equation*}

We finally choose $\ell=2\lambda$ and \begin{equation*}
    w:= \left\{
\begin{array}{l}
    \rho_{\lambda,\delta}e^{i\varphi_\lambda}\quad \text{if }\gq >0,\\
    \rho_{\lambda,\delta}e^{-i\varphi_\lambda}\quad \text{if }\gq <0\\
    1\quad\text{if }\gq =0.
    \end{array}
\right.
\end{equation*}

When $\gq\neq 0$, all conditions are fulfilled with the choice of $\delta$, except for the estimate \eqref{lemme  E(w) leq C |q|}. In view of the proof of Lemma 6 in~\cite{BetGrSa2}, it remains to deal with the potential energy. Since $\rho_{\lambda,\delta}$ is bounded, we use a Taylor expansion of $F$ near 1 and we get that

\begin{align*}
    F(\rho_{\lambda,\delta}^2) & = \big| F(1 -\delta-f_\lambda)\big|\\
    & =\big| F(1 -\delta-f_\lambda) - F(1) - F'(1)(-\delta-f_\lambda)\big|\\
    &\leq \dfrac{\Vert f'\Vert_{L^\ii([-3,3])}}{2}(\delta+f_\lambda)^2,
\end{align*}

then, reducing $\delta$ if necessary, we obtain

\begin{align*}
    E(w) & \leq \int_0^{2\lambda}\bigg( \dfrac{(f_\lambda ')^2}{8(1-\delta-f_\lambda)}+\big( 1-\delta-f_\lambda\big) \dfrac{(\varphi_\lambda')^2}{2} +C\Big( \dfrac{f_\lambda^2}{4}+\dfrac{\delta f_\lambda}{2}+\dfrac{\delta^2}{4}\Big)\bigg)\\
    & \leq \dfrac{1}{4\lambda^3}+\dfrac{1}{\lambda} + C\Big( \frac{1}{48\lambda}+\dfrac{\delta}{8}+\dfrac{\delta^2 \lambda}{2}\Big).
\end{align*}

We conclude as in Lemma 6 in~\cite{BetGrSa2}, by taking $\delta$ small enough and proportional to $\frac{1}{\lambda}$. 
\end{proof}

\begin{lem}\label{lemme compacité des suites pseudomin dans H^1 loc et inegalité faible}
Let $\gp$ satisfying~\eqref{hypothèse sur gp en fonction de q_*} and let $(u_n)$ be a pseudo-minimizing sequence. Then, there exist a subsequence $(u_{\sigma(n)})$ and a solution $v_c$ of \eqref{TWC} such that 

\begin{equation*}
    u_{\sigma(n)}\ntendf v_c\quad\text{in }H^1_{\loc}(\R).
\end{equation*}

Either $v_c$ is a constant map of modulus 1 or a non constant travelling wave that we relabel $\gv_c$.

Moreover, for any $A>0$,

\begin{equation}\label{lem: convergence sur [-A,A] de l'énergie} 
    \int_{-A}^A e(v_c)\leq \liminf_{n\rightarrow + \ii}\int_{-A}^{A} e(u_n)
\end{equation}

and if $v_c$ does not vanish on $\R$, then

\begin{equation}\label{lem: convergence sur [-A,A] du moment} 
    \int_{-A}^A (1-\rho_c^2)\varphi_c'\leq \lim_{n\rightarrow + \ii}\int_{-A}^{A}  (1-\rho_n^2)\varphi_n',
\end{equation}

where we have written $u_n=\rho_n e^{i\varphi_n}$ and $v_c =\rho_c e^{i\varphi_c}$.
\end{lem}

\begin{rem}
    If we suppose that the hypotheses of Theorem~\ref{existence et unicité des travelling wave si condition du zéro} holds, then we also know that the travelling wave~$\gv_c$ that we obtain is unique, up to a translation and a constant phase shift.
\end{rem}

\begin{rem}\label{remarque suivant le lemme compacité des suites pseudomin dans H^1 loc}
By letting $A$ tend to $+\ii$ in~\eqref{lem: convergence sur [-A,A] de l'énergie}, we can deduce that $E(v_c)\leq \Emin(\gp)$.
    In particular, if $\gp <\gq_*$, then by Lemma~\ref{lemme Emin est strictement croissante}, we get $E(v_c)\leq\Emin (\gp)<\Emin (\gq_*)=E(\gv_0)$ thus $v_c$ does not vanish. 
\end{rem}

\begin{rem}\label{remarque: si v_c(0)=0, alors c=0}
    One shall notice that if $v_c$ is a solution as in Lemma~\ref{lemme compacité des suites pseudomin dans H^1 loc et inegalité faible} that vanishes, then it must be the black soliton. Otherwise, if $c\neq 0$, assuming that it vanishes at $x_0\in\R$, and writing the first order differential equation satisfied by $\eta_c:=1-|v_c|^2$. This yields to $-\big(\eta_c'(x_0)\big)^2 = \mathcal{N}_c\big(\eta_c(x_0)\big)=c^2 >0$. However, since $\eta_c(x_0)=1$, then $x_0$ is a global maximum of the function $\eta_c$, thus $\eta_c'(x_0)=0$, which leads to a contradiction.
\end{rem}

\begin{proof}
Since $\big(E(u_n)\big)_n$ is bounded, the sequence $(u_n')_n$ is also bounded in $L^2(\R)$. By hypothesis \eqref{hypothèse de croissance sur F minorant intermediaire}, and the fact that $\big(E(u_n)\big)_n$ is bounded, we get that $(1-|u_n|^2)_n$ is bounded in $L^2(\R)$. Thus, $(1-|u_n|)_n$ and then $(u_n)_n$ are bounded in $L^2_{\loc}(\R)$. We then conclude that there exists a function $u$ such that, up to a subsequence,

\begin{equation*}
    u_n\ntendf u\quad\text{in }H^1_{\loc}(\R). 
\end{equation*}

By Rellich compactness theorem, the strong convergence holds in $L^{\ii}_\loc(\R)$, still up to a subsequence and in particular, \eqref{lem: convergence sur [-A,A] de l'énergie} holds. If the limit $u$ does not vanish on any interval $[-A,A]$ and then we can lift it, the limit $u$ also satisfies the inequality \eqref{lem: convergence sur [-A,A] du moment}. It remains to verify that $u$ is a solution of \eqref{TWC}. For that purpose, we consider $\xi\in\mathcal{C}^\ii_c(\R)$ such that 
\begin{equation}\label{preuve lemme int langle iu,xi' rangle=0}
    \int_\R \langle iu,\xi'\rangle_\C =0.
\end{equation}

We shall prove that

\begin{equation*}
    \int_{-A}^A e(u_{\sigma(n)}+t\xi)\geq \int_{-A}^A e(u_{\sigma(n)}) + O(t^2) + \underset{n\rightarrow +\ii}{o(1)}. 
\end{equation*}

The functions $u_n$ do not vanish, so we can compute the momentum of $u_n+t\xi$, noticing that for $R$ large enough $u_n(\pm R)+t\xi(\pm R)=u_n(\pm R)$ because $\xi$ is compactly supported. We obtain

\begin{align*}
    p(u_n+t\xi)=p(u_n) + t\int_\R \langle i u_n ,\xi'\rangle +O(t^2),
\end{align*}

where $O(t^2)$ does not depend on $n$. We use assumption \eqref{preuve lemme int langle iu,xi' rangle=0} and the convergence in $H^1_{\loc}(\R)$ to state that 

\begin{equation*}
    \int_\R \langle i u_n ,\xi'\rangle = \underset{n\rightarrow +\ii}{o(1)}.
\end{equation*}

Therefore

\begin{equation*}
    p(u_n+t\xi)=\gp +O(t^2) + \underset{n\rightarrow +\ii}{o(1)},
\end{equation*}

so that, setting $\gq_{n,t} = \gp - p(u_n+t\xi)$, we get

\begin{equation*}
    \gq_{n,t}=O(t^2) + \underset{n\rightarrow +\ii}{o(1)}.
\end{equation*}

Following the same construction as in Lemma 7 in~\cite{BetGrSa2}, we invoke Lemma~\ref{lemme equivalent du lemme 6 dans existance and travelling Begrasa 2} with $\gq = \gq_{n,t}$ and $\mu=\mu_{n,t}:=\inf\{\frac{1}{4},\frac{\nu_{n,t}}{2}\}$ where $\nu_{n,t}:=\sup\big\{ \big| 1-|u_n(x)|\big| x\notin [-A,A]\big\}$. This yields a positive number $l_{n,t}>1$ and a map $w_{n,t}=|w_{n,t}|e^{i\varphi_{n,t}}$ defined on $[0,l_{n,t}]$ such that 
\begin{equation*}
    w_{n,t}(0)=w_{n,t}(l_{n,t})\quad\text{and}\quad\big| 1-|w_{n,t}(0)|\big| \leq \mu_{n,t}.
\end{equation*}

Moreover, we have \begin{equation*}
    \gq_{n,t}=\dfrac{1}{2}\int_0^{l_{n,t}} (1-|w_{n,t}|^2)\varphi_{n,t}',
\end{equation*}

and 
\begin{equation*}
    E(w_{n,t})\leq C|\gq_{n,t}|=O(t^2) + \underset{n\rightarrow +\ii}{o(1)}.
\end{equation*}

In view of the mean value theorem, there exists some point $x_n\in [A,+\ii)$ such that $|u_n(x_n)|=|w_{n,t}(0)|$. Multiplying possibly $w_{n,t}$ by a complex number of modulus one, we can assume that $u_n(x_n)=w_{n,t}(0)$. We define a comparison map as follows,

\begin{equation*}
    v_{n,t}(x):=\left\{
\begin{array}{l}
    u_n(x)+t\xi(x)\quad\text{if }x<x_n,\\
    w_{n,t}(x-x_n)\quad \text{if }x\in [x_n,x_n+l_{n,t}],\\
    u_n(x-l_{n,t})+t\xi(x-l_{n,t})\quad\text{if }x\geq x_n + l_{n,t}.
    \end{array}
\right.
\end{equation*}

We verify that $E(v_{n,t})=E(u_n+t\xi)+E(w_{n,t})$ and $p(v_{n,t})=p(u_n+t\xi)+p(w_{n,t})=\gp$, therefore

\begin{equation*}
    \Emin(\gp)\leq E(v_{n,t}).
\end{equation*}

Since $(u_n)$ is a pseudo-minimizing sequence, we have \begin{equation*}
    E(u_n)=\Emin(\gp) + \underset{n\rightarrow +\ii}{o(1)},
\end{equation*}

whereas, since $\xi$ has compact support, 

\begin{equation*}
    E(u_n+t\xi)-E(u_n) = \int_{-A}^A \Big(e(u_n+t\xi)-e(u_n)\Big).
\end{equation*}

We infer from the previous estimates that

\begin{align*}
    \int_{-A}^A \Big(e(u_{\sigma(n)}+t\xi)-e(u_{\sigma(n)})\Big) & = E(u_{\sigma(n)}+t\xi)-E(u_{\sigma(n)})\\
    & = E(v_{\sigma(n),t})-E(w_{\sigma(n),t}) - \big( \Emin (\gp) + \underset{n\rightarrow +\ii}{o(1)}\big)\\
    & \geq \Emin (\gp) + O(t^2) + \underset{n\rightarrow +\ii}{o(1)} - \big( \Emin (\gp) +\underset{n\rightarrow +\ii}{o(1)}\big) = O(t^2) +\underset{n\rightarrow +\ii}{o(1)},
\end{align*}

so that the claim is proved. To conclude, we expand the integral in the claim so that 

\begin{align*}
    \int_{-A}^A tu_n' \xi' + \dfrac{1}{2}\big( F(|u_n+t\xi|^2)-F(|u_n|^2)\big)  & \geq O(t^2) + \underset{n\rightarrow +\ii}{o(1)}.
\end{align*}

Letting $n$ tend to $+\ii$, this yields

\begin{align*}
    \int_{-A}^A tu' \xi' + \dfrac{1}{2}\big( F(|u+t\xi|^2)-F(|u|^2)\big)  & \geq O(t^2).
\end{align*}

By writing the first order Taylor expansion of $F$ in $|u|^2$, we obtain
\begin{align*}
    \int_{-A}^A t u' \xi' - t u \xi f(|u|^2)  & \geq O(t^2),
\end{align*}

so that

\begin{align*}
    -t\int_{-A}^A u \xi'' + u \xi f(|u|^2) & \geq O(t^2).\\
\end{align*}

Now letting $t$ tend to $0^+$ and $0^-$, we deduce that
\begin{equation*}
    \int_{-A}^A u'' \xi + u \xi f(|u|^2)=0 .
\end{equation*}

Since $\xi$ is any arbitrary function with compact support satisfying \eqref{preuve lemme int langle iu,xi' rangle=0}, this implies the existence of a constant $c\in\R$ such that $u$ solves \eqref{TWC}.
\end{proof}

Let $\gp$ satisfy~\eqref{hypothèse sur gp en fonction de q_*} and $(u_n)$ be a pseudo-minimizing sequence as in \eqref{etre une suite pseudo minimisante}. We set for $n$ large enough, $$\delta_n:= 1-\dfrac{E(u_n)}{c_s |p(u_n)|},$$

thus,

\begin{equation*}
    \delta_n\ntend \delta_\mathfrak{p}:=1-\dfrac{E_{\min}(\mathfrak{p})}{c_s\mathfrak{p}}.
\end{equation*}

One crucial remark is that for $\gp >0$, $\Emin (\gp) <  c_s\gp$ by Corollary \ref{corollaire Emin is striclty sub additive}. We infer that $\delta_{\gp}>0$, so that for $n$ large enough, $\delta_n \geq \frac{\delta_\mathfrak{p}}{2}$. We are almost in position to prove Theorem \ref{theoreme de continuité de l'energie dans le cas d'une non linéarité générale}. When $n$ is sufficiently large, we can apply Lemma~\ref{lemme controle de v pres et loin de 1} with $E=\frac{c_s \gp+\Emin (\gp)}{2}$ and $\delta_0=\frac{\delta_{\gp}}{4}$. Therefore, either for all $n$,

\begin{equation*}
    \big| 1 - |u_n(x)|\big| < \dfrac{\delta_{\gp}}{4}\quad\text{for all }x\in\R,
\end{equation*}

or there exists an integer $l_{\gp}$ only depending on $\gp$ such that there exist $l_n$ points $x_1^n ,..., x_{l_n}^n$ with $l_n \leq l_{\gp}$ such that, the second property of Lemma~\ref{lemme controle de v pres et loin de 1} holds. For $n$ large enough, we also know that $1-\frac{E(u_n)}{c_s \gp}>0$. Because of \eqref{corollaire: 1-|v(x_delta)| geq delta} in Corollary~\ref{corollaire du lemme sqrt(f rho^2)|varphi'| leq lambda e(v)/rho} with $\lambda=\frac{c_s^2}{4}$, there exists $x\in\R$ such that $1-|u_n(x)| \geq \frac{\delta_{\gp}}{4}$. By exhaustion of the first case in Lemma~\ref{lemme controle de v pres et loin de 1}, we can only have

\begin{equation*}
    \big| 1 - |u_n(x^n_i)|\big|\geq\dfrac{\delta_\mathfrak{p}}{4}, \quad\forall i\in \lbrace 1 ,..., l_n\rbrace,
\end{equation*}

and

\begin{equation*}
    \big| 1 - |u_n(x)|\big|\leq\dfrac{\delta_\mathfrak{p}}{4}, \quad\forall x \in\R\setminus \bigcup_{i=1}^{l_n} [x_i^n-1, x_i^n + 1 ].
\end{equation*}

Passing possibly to a subsequence, we can assume that the number $l_n$ does not depend on $n$, and set $l=l_n$. A standard compactness argument shows, that passing again possibly to a further subsequence, and relabelling possibly the points $x_i^n$, we may find some integer $1\leq \widetilde{l}\leq l$ and $R>0$ such that 

\begin{equation*}\label{x_i^n et x_j^n s'éloignent quand n grandit}
    |x_i^n - x_j^n|\ntend +\ii,\quad\forall 1\leq i\neq j\leq \widetilde{l},
\end{equation*}

and 

\begin{equation*}
    x_i^n \in\bigcup_{j=1}^{\widetilde{l}} \mathcal{B}(x_j^n,R),\quad\forall \widetilde{l} < i\leq l.
\end{equation*}

We deduce that
\begin{equation*}
    \big| 1 - |u_n(x)|\big|\leq\dfrac{\delta_\mathfrak{p}}{4}, \quad\forall x \in\R\setminus \bigcup_{i=1}^{\widetilde{l}} \mathcal{B}(x_i^n,R+1).
\end{equation*}

By the inequality \eqref{cor: 1ere inegalité} in Corollary \ref{corollaire du lemme sqrt(f rho^2)|varphi'| leq lambda e(v)/rho}, we obtain

\begin{equation}\label{inégalité avec C_0 qui permet de faire le step 3}
    \dfrac{1}{2}\big| (|u_n|^2-1)\varphi_n'\big|\leq \dfrac{1}{c_s}\dfrac{e(u_n)}{1-\frac{\delta_\mathfrak{p}}{4}}=\dfrac{e(u_n)}{C_0},\quad \forall x \in\R\setminus \bigcup_{i=1}^{\widetilde{l}} \mathcal{B}(x_i^n,R+1).
\end{equation}

where $C_0:=c_s (1-\frac{\delta_\mathfrak{p}}{4})$ by \eqref{hypothèse de croissance sur F minorant intermediaire}.

\begin{lem}\label{inf |u_n| geq alpha_0}
    For $\gp$ satistying~\eqref{hypothèse sur gp en fonction de q_*}, there exists $\alpha_0>0$ independent of $n$ such that 

\begin{equation}\label{lemme: inf|u_n| geq alpha_0}
    \inf_\R |u_n|\geq \alpha_0 .
\end{equation}
\end{lem}

\begin{proof}
We suppose by contradiction that there exists $(a_n)$ such that 

\begin{equation}\label{u_n(a_n)ntend 0}
    u_n(a_n)\ntend 0.
\end{equation}

We apply Lemma~\ref{lemme compacité des suites pseudomin dans H^1 loc et inegalité faible} to $\big( u_n(.+a_n)\big)$, we label the limit function $v_c$ and we consider the two following cases. Combining~\eqref{u_n(a_n)ntend 0} and the Rellich compactness theorem, we have $v_c(0)=0$, thus there exist $\theta,\widetilde{x}\in\R$ such that $v_c=e^{i\theta}\gv_0(.+\widetilde{x})$ by Remark~\ref{remarque: si v_c(0)=0, alors c=0}. If $\gp <\gq_*$, no matter the value of $\gq_*$, we obtain a contradiction with Remark~\ref{remarque suivant le lemme compacité des suites pseudomin dans H^1 loc}. Now, if $\gp=\gq_*$, we have
\begin{equation*}
    E(u_n)\ntend E(v_0),
\end{equation*}
and this will provide strong convergences for $(u_n)$. We will then infer that 
\begin{equation}\label{lemme preuve convergence p(u_n) ntend [p](v_0)}
    p(u_n)=[p](u_n)\ntend [p](\gv_0)=\frac{\pi}{2},
\end{equation}
which contradicts the fact that $\gp\neq \frac{\pi}{2}\mod\pi$.~\footnote{One shall keep in mind the following argument, since it will be used to prove the convergences in the concentration-compactness theorem.}\\
More precisely, for $\mu >0$, there exist $A >0$ and $N\in\N$ such that, if $n\geq N$, 

\begin{equation*}
    \int_{-A+a_n}^{A+a_n} e(u_n) \geq E(\mathfrak{v}_0)-\mu 
\end{equation*}

and

\begin{equation*}
     E(\mathfrak{v}_0)-\mu -\dfrac{1}{2}\int_{-A}^A F\big( |u_n(.+a_n)|^2\big)\leq \dfrac{1}{2}\Vert  u'_n(.+a_n)\Vert^2_{L^2(\R)} \leq E(u_n)-\dfrac{1}{2}\int_{-A}^A F\big( |u_n(.+a_n)|^2\big).
\end{equation*}

By passing to the limit $n\rightarrow +\ii$, this yields

\begin{equation}\label{inégalité après avoir pris limite inf et sup}
    E(\mathfrak{v}_0)-\mu-\dfrac{1}{2}\int_{-A+\widetilde{x}}^{A+\widetilde{x}}F(|\mathfrak{v}_0|^2)\leq \dfrac{1}{2}\liminf_{n\rightarrow +\ii}\Vert u'_n(.+a_n)\Vert^2_{L^2(\R)}\leq E(\mathfrak{v}_0)-\dfrac{1}{2}\int_{-A+\widetilde{x}}^{A+\widetilde{x}}F(|\mathfrak{v}_0|^2) ,
\end{equation}

and the same inequality holds with the limsup. Furthermore, for $A$ large enough, one gets,

\begin{equation*}
    \dfrac{1}{2}\int_\R |\mathfrak{v}'_0|^2 -\mu\leq E(\mathfrak{v}_0)-\dfrac{1}{2}\int_{-A+\widetilde{x}}^{A+\widetilde{x}}F(|\mathfrak{v}_0|^2)\leq \dfrac{1}{2}\int_\R |\mathfrak{v}'_0|^2+\mu .
\end{equation*}

Introducing both these inequalities into \eqref{inégalité après avoir pris limite inf et sup}, we have

\begin{equation*}
    \dfrac{1}{2}\int_\R |\mathfrak{v}'_0|^2 -2\mu\leq \dfrac{1}{2}\liminf_{n\rightarrow +\ii}\Vert u'_n(.+a_n)\Vert^2_{L^2(\R)}\leq \dfrac{1}{2}\limsup_{n\rightarrow +\ii}\Vert u'_n(.+a_n)\Vert^2_{L^2(\R)}\leq  \dfrac{1}{2}\int_\R |\mathfrak{v}'_0|^2+\mu .
\end{equation*}

Since this is true for any $\mu >0$, we finally get 

\begin{equation*}
    \lim_{n\rightarrow +\ii}\Vert u'_n\Vert_{L^2} =\Vert \mathfrak{v}'_0\Vert_{L^2}.
\end{equation*}

Therefore $u'_n(.+a_n)\ntend \mathfrak{v}'_0(.+\widetilde{x})$ strongly in $L^2(\R)$. Since $u_n(.+a_n)\ntendf e^{i\theta} \mathfrak{v}_0(.+\widetilde{x})$ weakly in $H^1_{\loc}(\R)$, up to a subsequence we also know that $u_n(.+a_n)\ntend e^{i\theta} \mathfrak{v}_0(.+\widetilde{x})$ 
uniformly on any compact subset of $\R$. Finally, we additionally have \begin{equation}\label{F(w_n) tend vers F(v_0) dans L^1}
    F(|u_n(.+a_n)|^2)\ntend F(|\mathfrak{v}_0(.+\widetilde{x})|^2)\quad\text{in }L^1(\R),
\end{equation}

since $E(u_n)\ntend E(\mathfrak{v}_0)$. Still by the local uniform convergence,

\begin{equation*}
    1-|u_n(.+a_n)|^2 \ntend 1- |\gv_0(.+\widetilde{x})|^2\quad a.e.,
\end{equation*}

and we have, using assumption \eqref{hypothèse de croissance sur F minorant intermediaire},

\begin{align*}
    \big| 1-|u_n(.+a_n)|^2\big|^2 \leq \dfrac{1}{\lambda}F\big(|u_n(.+a_n)|^2\big).
\end{align*}

By \eqref{F(w_n) tend vers F(v_0) dans L^1}, we can assume, up to a subsequence, that there exists $h\in L^1(\R)$ such that for all $n\in\N$, 

$$F\big(|u_n(.+a_n)|^2\big) \leq h\quad\text{a.e.},$$

so that we get the upper bound

\begin{align*}
    \big| 1-|u_n(.+a_n)|^2\big|^2  & \leq \dfrac{1}{\lambda}h.
\end{align*}

By the dominated convergence theorem, we obtain

\begin{equation*}
    1-|u_n(.+a_n)|^2\ntend 1- |\gv_0(.+\widetilde{x})|^2\quad\text{in }L^2(\R),
\end{equation*}

and this yields in particular
\begin{equation*}        d_A\big(u_n(.+a_n),\gv_0(.+\widetilde{x})\big)\ntend 0,
\end{equation*}

so that, by Lemma~\ref{lemme appendix untwisted momentum continu pour d_A}, we obtain~\eqref{lemme preuve convergence p(u_n) ntend [p](v_0)} and this concludes the proof.

\end{proof}

Now we can prove Theorem \ref{theoreme de continuité de l'energie dans le cas d'une non linéarité générale}.

\begin{proof}[Proof of Theorem~\ref{theoreme de continuité de l'energie dans le cas d'une non linéarité générale}]
\textit{Step 1}. For any $1 \leq i\leq\widetilde{l}$, there exists $c_i \in (0,c_s)$, such that 

\begin{equation*}
    u_n(.+x^n_i)\ntendf \mathfrak{v}_{c_i}\quad\text{in }H^1_{\loc}(\R).
\end{equation*}

Indeed, applying Lemma \ref{lemme compacité des suites pseudomin dans H^1 loc et inegalité faible} to the sequences $\big( u_n(. + x_i^n)\big)_n$ yields the existence of $v_{c_i}$, a limiting solution to $(TW_{c_i})$ with $c_i\in\R$. Furthermore we have 

\begin{equation*}
    |u_n(x_i^n)|\leq 1-\dfrac{\delta_\gp}{4},
\end{equation*}

we deduce from the uniform convergence on every compact set of $\R$ that \begin{equation*}
    |v_{c_i}(0)|\leq 1-\dfrac{\delta_\gp}{4}.
\end{equation*}

We infer that $v_{c_i}$ is not a constant map with moreover $c_i\in [0,c_s)$. In order to be consistent with the statement of Theorem~\ref{theoreme de continuité de l'energie dans le cas d'une non linéarité générale}, we relabel this travelling wave as $\gv_{c_i}$.

Furthermore, we use Lemma~\ref{inf |u_n| geq alpha_0} and the local uniform convergence to show that $c_i\neq 0$. Indeed, if by contradiction $c_i=0$, then by uniqueness, $\gv_{c_i}$ is a translation of the black soliton (up to a phase change). The travelling wave $\gv_{c_i}$ thus vanishes at some point $x_0$. For $n$ large enough, we infer by local uniform convergence, that $|u_n(x_0+x_i^n) - \gv_{c_i}(x_0)|\leq \frac{\alpha_0}{2}$, with $\alpha_0$ as in~\eqref{lemme: inf|u_n| geq alpha_0}. Therefore

\begin{align*}
    0=|\gv_{c_i}(x_0)|  &\geq |u_n(x_0+x_i^n)|- |\gv_{c_i}(x_0)-u_n(x_0+x_i^n)|\geq\dfrac{\alpha_0}{2}.
\end{align*}

This brings a contradiction, therefore we conclude that $c_i\neq 0$.

\textit{Step 2}. Given any number $\mu >0$, there exist numbers $A_\mu >0$ and $n_\mu \in\N^*$ such that

\begin{equation*}
    n\geq n_\mu \Longrightarrow \left\{
\begin{array}{l}
    \displaystyle\int_{\bigcup_{i=1}^{\widetilde{l}} \mathcal{B}(x_i^n, A_\mu)} e(u_n) \geq \sum_{i=1}^{\widetilde{l}} E(\mathfrak{v}_{c_i})-\mu \\
    \Big|\dfrac{1}{2}\displaystyle\int_{\bigcup_{i=1}^{\widetilde{l}} \mathcal{B}(x_i^n, A_\mu)}(\rho_n^2 - 1 )\varphi_n' - \sum_{i=1}^{\widetilde{l}}\gp_i\Big|\leq \mu .
    \end{array}
\right.
\end{equation*}

where $\gp_i := p(\mathfrak{v}_{c_i})\neq 0$. In view of Step 1, we only need to take $A_\mu > R+1$ \footnote{We will use this condition in Step 3. It helps us to control the momentum in terms of the energy on the area where $|u_n|$ is concentrated near $1$.} such that, for any $1\leq i\leq \widetilde{l}$, we have

\begin{equation*}
    E(\mathfrak{v}_{c_i})\leq \int_{-A_\mu+\widetilde{x}_i}^{A_\mu+\widetilde{x}_i} e(\mathfrak{v}_{c_i}) +\dfrac{\mu}{2\widetilde{l}}\quad\text{and}\quad \Big| \dfrac{1}{2}\int_{-A_\mu+\widetilde{x}_i}^{A_\mu+\widetilde{x}_i} (\rho_{c_i}^2-1)\varphi_{c_i}'-\gp_i\Big|\leq\dfrac{\mu}{2\widetilde{l}}
\end{equation*}

with $R>0$ exhibited earlier and such that \eqref{inégalité avec C_0 qui permet de faire le step 3} holds. We deduce that, for $n$ large enough,

\begin{equation*}
    \int_{-A_\mu+x_i^n}^{A_\mu+x_i^n} e(u_n) \geq E(\mathfrak{v}_{c_i})-\dfrac{\mu}{\widetilde{l}}\quad\text{and}\quad
    \Big|\dfrac{1}{2}\displaystyle\int_{-A_\mu+x_i^n}^{A_\mu + x_i^n} (\rho_{n}^2 - 1 )\varphi_{n}' - \gp_i\Big|\leq \dfrac{\mu}{\widetilde{l}}\ .
\end{equation*}

Step 2 then follows from summing.

\textit{Step 3}. We have 

\begin{equation*}
    \Big| \dfrac{1}{2}\int_{\R\setminus \bigcup_{i=1}^{\widetilde{l}} \mathcal{B}(x_i^n, A_\mu)}(\rho_n^2-1)\varphi_n'\Big|\leq \dfrac{1}{C_0}\int_{\R\setminus \bigcup_{i=1}^{\widetilde{l}} \mathcal{B}(x_i^n, A_\mu)} e(u_n).
\end{equation*}

Since we took $A_\mu > R+1$, this is just a consequence of integrating \eqref{inégalité avec C_0 qui permet de faire le step 3}. Moreover, since $\big(E(u_n)\big)_n$ is bounded, passing possibly to a further sequence, we can suppose that there exist $E_\mu$ and $\mathfrak{p}_\mu$ such that

\begin{equation*}
    \int_{\R\setminus \bigcup_{i=1}^{\widetilde{l}} \mathcal{B}(x_i^n, A_\mu)}(\rho_n^2-1)\varphi_n'\rightarrow \gp_\mu\text{ and }\int_{\R\setminus \bigcup_{i=1}^{\widetilde{l}} \mathcal{B}(x_i^n, A_\mu)} e(u_n)\rightarrow E_\mu .
\end{equation*}

Going back to Step 2, and letting $n\rightarrow +\ii$, we are led to 

\begin{equation*}
    \Emin(\gp)\geq \displaystyle\sum_{i=1}^{\widetilde{l}}\Emin(\gp_i)  + E_\mu - \mu \text{ and } \Big| \gp - \displaystyle\sum_{i=1}^{\widetilde{l}} \gp_i - \gp_{\mu} \Big| \leq \mu,   
\end{equation*}

whereas Step 3 yields

\begin{equation*}
    C_0 |\gp_\mu| \leq E_\mu .
\end{equation*}

We may assume that, up to a subsequence $(\mu_m)$ tending to 0, we have 

\begin{equation*}
    \gp_{\mu_m}\underset{m\rightarrow +\ii}{\longrightarrow}\widetilde{\gp}\text{ and }E_{\mu_m}\underset{m\rightarrow +\ii}{\longrightarrow} \widetilde{E},
\end{equation*}

which finally leads to

\begin{equation}\label{sum p_i +widetilde = p ou pi/2 + kpi}
    \sum_{i=1}^{\widetilde{l}}\gp_i + \widetilde{\gp} = 
    \gp,
\end{equation}

\begin{equation}\label{Emin(p) geq sum Emin(p_i) + Etilde}
   \Emin(\gp)\geq \displaystyle\sum_{i=1}^{\widetilde{l}}\Emin(\gp_i)  + \widetilde{E},
\end{equation}

and 

\begin{equation}\label{C_0 p leq Etilde}
    C_0 |\widetilde{\mathfrak{p}}|\leq \widetilde{E}.
\end{equation}

\textit{Step 4}.
We prove that $\widetilde{E}=\widetilde{\gp}=0$ and $\widetilde{l}=1$. Observe that

\begin{equation}\label{Emin(p)/p leq C_0}
        \dfrac{\Emin(\gp)}{\gp}=\mathfrak{c}_s (1-\delta_\mathfrak{p}) <\mathfrak{c}_s (1-\dfrac{\delta_\mathfrak{p}}{4})= C_0.
\end{equation}

By contradiction, we suppose first that $\widetilde{\gp}\neq 0$. Then using successively the evenness of $\Emin$, then~\eqref{Emin(p) geq sum Emin(p_i) + Etilde}, \eqref{C_0 p leq Etilde}, \eqref{Emin(p)/p leq C_0}, the concavity, the monotonicity of $\Emin$ and \eqref{sum p_i +widetilde = p ou pi/2 + kpi}, we write

\begin{align*}    \Emin(\gp)\geq\sum_{i=1}^{\widetilde{l}}\Emin(|\gp_i|) + \widetilde{E} > \sum_{i=1}^{\widetilde{l}}\Emin(|\gp_i|) +\dfrac{\Emin(\gp)}{\gp}|\widetilde{\gp}| \geq \dfrac{\Emin(\gp)}{\gp}\Big(\sum_{i=1}^{\widetilde{l}}|\gp_i| + |\widetilde{\gp}|\Big)\geq \Emin(\gp),
\end{align*}

so that $\widetilde{\gp}=0$. By~\eqref{sum p_i +widetilde = p ou pi/2 + kpi}, we have

\begin{equation}\label{sum emin(p_i) + widetilde E less ...}    \sum_{i=1}^{\widetilde{l}}\Emin(|\gp_i|) + \widetilde{E} \leq \Emin(\gp)=\Emin\Big(\sum_{i=1}^{\widetilde{l}}\gp_i\Big) \leq \sum_{i=1}^{\widetilde{l}}\Emin(|\gp_i|),
\end{equation}

then $\widetilde{E}\leq 0$, hence $\widetilde{E}=0$. Finally, if $\widetilde{l}\geq 2$, by strict subadditivity, we should obtain $$\Emin\Big(\sum_{i=1}^{\widetilde{l}}\gp_i\Big) < \sum_{i=1}^{\widetilde{l}}\Emin(|\gp_i|),$$

which happens to be an equality in view of~\eqref{sum emin(p_i) + widetilde E less ...}. We conclude that $\widetilde{E}=\widetilde{\gp}=0$ and $\widetilde{l}=1$.

\textit{Step 5}. Conclusion. By Step 4, we know that there exist $c_{i_0}\in (0,c_s)$ and a sequence $(x^n_{i_0})$ such that, 
\begin{equation*}
    u_n(.+x^n_{i_0})\underset{n\rightarrow +\ii}{\rightharpoonup} \mathfrak{v}_{c_{i_0}}\quad\text{in } H^1_{\loc}(\R)\quad\text{and}\quad p(\mathfrak{v}_{c_{i_0}})=\gp.
\end{equation*}

Since $\widetilde{l}=1$, we relabel the parameters $c_{i_0}$ and $x^n_{i_0}$ as $c$ and $x_n$. We have shown that $p(\gv_c)=\gp$. This provides, in addition to Remark~\ref{remarque suivant le lemme compacité des suites pseudomin dans H^1 loc}, the estimate

\begin{equation*}
    \Emin(\gp)\leq E(\gv_c).
\end{equation*}

Thus the previous inequality is in fact an equality and this is why

\begin{equation*}\label{gv_c in S_gp}
\gv_c\in\mathcal{S}_\gp.
\end{equation*}

Now, replacing $\gv_0$ by $\gv_c$ in the energetic argument in the proof of Lemma~\ref{inf |u_n| geq alpha_0}, we recover the desired convergences and this concludes the proof of Theorem~\ref{theoreme de continuité de l'energie dans le cas d'une non linéarité générale}. 
\end{proof}

We conclude that there exist minimizers of the energy when the momentum is fixed. These are necessarily travelling waves.

\begin{rem}
    Although we obtain, up to an extraction of the pseudo-minimizing sequence, the convergence to a travelling wave $\gv_c$, we have no information regarding the speed $c$. This is one of the drawbacks of such a compactness method. As a consequence, we will see that we are not willing to prove the orbital stability of a travelling wave when its speed is fixed in advance. By contrast, this theorem provides the existence of travelling waves with given momentum.
\end{rem}

\section{Orbital stability}\label{section stabilité orbitale}
In this section, we finish the proof of Theorem~\ref{théorème principal}. We recall that

\begin{equation*}
    \energyset:=\big\lbrace w\in L^\ii(\R)\big| w'\in L^2(\R), F(|w|^2)\in L^1(\R)\big\rbrace .
\end{equation*}

For $A >0$, we consider the distance

\begin{equation*}
    d_A(v_1,v_2):=\Vert v_1-v_2\Vert_{L^{\ii}([-A,A])} +\Vert v_1'-v_2'\Vert_{L^2}+\Vert |v_1|^2-|v_2|^2\Vert_{L^2},
\end{equation*}

and we prove the orbital stability of $\mathcal{S}_\gp$ when $\gp$ satisfies \eqref{hypothèse sur gp en fonction de q_*}.

\begin{proof}
We assume by contradiction that we may find $\varepsilon_0 >0$ such that for all $n\in\N$, there exist sequences $(\delta_n),(t_n)\in\R^{\N}$ and $(\Psi_0^n)_n\in \energyset^\N$ such that $(\delta_n)$ tends to $0$, 

\begin{equation}\label{preuve du theoreme de stabilité orbitale, distance de psi_0 à v_c inférieure à delta_n}
    d_A(\Psi_0^n,\mathcal{S}_\gp)\leq\delta_n,
\end{equation}

and for all $a,\theta\in\R$,

\begin{equation}\label{preuve du theoreme de stabilité orbitale, d(psi(t),v_c(.-a)) > eps_0}
    d_A\Big(e^{i\theta}\Psi^n(t_n,.-a), \mathcal{S}_\gp\Big)>\varepsilon_0,
\end{equation}

where $\Psi^n$ denotes the solution of \eqref{NLS} corresponding to the initial datum $\Psi_0^n$. In particular, since we have \eqref{preuve du theoreme de stabilité orbitale, distance de psi_0 à v_c inférieure à delta_n}, for any $n\in\N$, there exists $v_n\in \mathcal{S}_\gp$ such that 

\begin{equation}\label{preuve du theoreme stab orb, distance de psi_0 a v_n inf a delta_n}
    d_A(\Psi_0^n,v_n)\leq 2\delta_n .
\end{equation}

Set $w_n:= \Psi^n(t_n)$. By conservation of the energy and of the momentum on the energy set $\energyset$, for all $n\in\N$, we have

\begin{equation*}
    E(w_n)=E(\Psi_0^n)\quad\text{and}\quad p(w_n)=p(\Psi_0^n).
\end{equation*}

In addition, we claim that 

\begin{equation}\label{theoreme preuve w_n est une suite pseudo min}
    E(\Psi_0^n)\ntend\Emin(\gp)\quad\text{and}\quad p(\Psi_0^n)\ntend\gp,
\end{equation}
so that $(w_n)_n$ is a pseudo-minimizing sequence.
By definition, $(v_n)$ is a pseudo-minimizing sequence, therefore by the concentration-compactness Theorem \ref{theoreme de continuité de l'energie dans le cas d'une non linéarité générale}, there exist $(a_n),\theta$ in $\R$ and $\gv_c\in \mathcal{S}_\gp$ such that 

\begin{equation}\label{w_n tend vers v_c dans C^0_loc}
    v_n(.+a_n)\ntend e^{i\theta} \gv_c\quad\text{in }\mathcal{C}^0_{\loc}(\R),
\end{equation}

\begin{equation}\label{w_n' tend vers v_c' dans L^2}
    v'_n(.+a_n)\ntend e^{i\theta} \gv_c'\quad\text{in }L^2(\R),
\end{equation}

and 

\begin{equation}\label{F(w_n) tend vers F(v_c) dans L^1}
    F\big(|v_n(.+a_n)|^2\big)\ntend F( |\gv_c|^2)\quad\text{in }L^1(\R).
\end{equation}

We are going to deduce \begin{equation}\label{d_A(pis_0^n,v) ntend 0}
d_A\big(\Psi_0^n(.+a_n),e^{i\theta}\gv_c\big)\ntend 0 .
\end{equation}

We have $d_A\big(\Psi_0^n(.+a_n),e^{i\theta}\gv_c\big)\leq d_A\big(\Psi_0^n(.+a_n),v_n(.+a_n)\big)+d_A\big(v_n(.+a_n),e^{i\theta}\gv_c\big)$.

First, there exists $C>0$ independent of $n$, such that 

\begin{equation}\label{d_A ( psi(.+a_n,...) leq C d_A ( psi(.,...)}
    d_A\big(\Psi_0^n(.+a_n),v_n(.+a_n)\big)\leq C d_A(\Psi_0^n,v_n).
\end{equation}

Indeed, we use the fact that $(a_n)$ is necessarily bounded by a positive number $M$. Otherwise, $\big(|v_n(x+a_n)|\big)$ would tend to 1 as $n$ tends to $+\ii$ for $x$ in any compact set whereas it tends to $|\gv_c(x)|$ that cannot be constant. Therefore, 

\begin{equation*}
    d_A\big(\Psi_0^n(.+a_n),v_n(.+a_n)\big)\leq d_{A+M}(\Psi_0^n,v_n),
\end{equation*}

and

\begin{align*}
    \Vert \Psi_0^n-v_n\Vert_{L^{\ii}([-A-M,A+M])}& \leq \Vert\Psi_0^n-v_n\Vert_{L^{\ii}([-A,A])}+\sqrt{A+M}\normLdeux{(\Psi_0^n)'-v_n'}\\
    &\leq C(A,M)d_A(\Psi_0^n,v_n).
\end{align*}

It remains to prove that

\begin{equation}\label{convergence de d_A vers 0}
d_A\big(v_n(.+a_n),e^{i\theta}\gv_c\big)\ntend 0.
\end{equation}

In view of the convergences~\eqref{w_n tend vers v_c dans C^0_loc}, \eqref{w_n' tend vers v_c' dans L^2} and~\eqref{F(w_n) tend vers F(v_c) dans L^1}, we only verify that

$$\Vert |v_n(.+a_n)|^2-|\gv_c|^2\Vert_{L^2}\ntend 0.$$

We can check this last statement using assumption \eqref{hypothèse de croissance sur F minorant intermediaire} and the dominated convergence theorem. Combining the previous assertion~\eqref{convergence de d_A vers 0} with~\eqref{d_A ( psi(.+a_n,...) leq C d_A ( psi(.,...)}, we obtain~\eqref{d_A(pis_0^n,v) ntend 0}. By Lemma~\ref{lemme convergence pour d_A implique la convergence dans H^1 de 1-rho_n^2 vers 1-rho^2}, we infer that 

\begin{equation*}
    p(\Psi_0^n)\ntend p(\gv_c).
\end{equation*}

In order to see that $E(\Psi_0^n)\ntend E(\gv_c)$, we first notice that $\Vert (\Psi_0^n)'\Vert_{L^2}\ntend \Vert \gv_c'\Vert_{L^2}$ by \eqref{convergence de d_A vers 0}. It remains to prove that

\begin{equation*}
    \Vert F(|\Psi_0^n|^2)\Vert_1\ntend\Vert F(|\gv_c|^2)\Vert_1.
\end{equation*}

We invoke Lemma \ref{lemme convergence pour d_A implique la convergence dans H^1 de 1-rho_n^2 vers 1-rho^2} to infer that

\begin{equation*}
    1-|\Psi_0^n|^2\ntend 1-|\gv_c|^2\quad\text{in }H^1(\R).
\end{equation*}

The proof of Lemma \ref{lemme convergence de 1-rho_n dans H^1 implique convergence de F(rho_n^2)} can then be adapted by replacing $\rho_n$ (resp. $\rho$) by $\rho_n^2$ (resp. $\rho^2)$ and it provides the convergence 

\begin{equation*}
    F(|\Psi_0^n|^2)\ntend F(|\gv_c|^2)\quad\text{in }L^1(\R).
\end{equation*}

Thus, we have showed that $(\Psi_0^n)$ and then $(w_n)$ are pseudo-minimizing sequences. Reasoning as before, we exhibit a sequence $(b_n)$, a real number $\widetilde{\theta}$ and a function $w\in\mathcal{S}_\gp$ that satisfy the statement of the concentration-compactness theorem. It provides 

\begin{equation*}
d_A\big( w_n(.+b_n),e^{i\widetilde{\theta}}w)\ntend 0,
\end{equation*}

thus

\begin{equation*}
    d_A\big( (e^{-i\widetilde{\theta}}w_n(.+b_n),\mathcal{S}_\gp\big)\ntend 0,
\end{equation*}

which contradicts \eqref{preuve du theoreme de stabilité orbitale, d(psi(t),v_c(.-a)) > eps_0} and then concludes the proof.

\end{proof}

\section{Numerical simulations}\label{section: numerical simulations}

In this section, we display some examples of the curves of $E(\gv_c),p(\gv_c)$ with respect to $c\in [0,c_s)$ according to the exact formulae\footnote{See~\cite{Chiron7} for more details.}:

\begin{equation}
    E(\gv_c)=2\int_0^{\xi_c}F(1-\xi)\dfrac{d\xi}{\sqrt{-\mathcal{N}_c(\xi)}}\quad\text{and}\quad p(\gv_c)=\dfrac{c}{2}\int_0^{\xi_c}\dfrac{\xi^2}{1-\xi}\dfrac{d\xi}{\sqrt{-\mathcal{N}_c(\xi)}}.
\end{equation}

For every example, we also plot the so-called energy/momentum diagram. This will give us a plain and precise idea of the localisation of $\gq_*$ in different cases. When the assumptions of Theorem~\ref{existence et unicité des travelling wave si condition du zéro}, are achieved for any $c$ i.e. when there is a unique travelling wave for each speed $c\in (0,c_s)$, we expect $\gq_*\geq \frac{\pi}{2}$ and we will show a partial result in that way in Appendix~\ref{subsection localisation de q_*}. Besides, we will also exhibit two examples in this direction. On the contrary, we will also exhibit a nonlinearity where $\gq_*$ seems to be less than $\frac{\pi}{2}$.

Comparatively to the next examples, we first exhibit the case of the Gross-Pitaevskii nonlinearity $f(\rho)=1-\rho$.

\begin{center}
    \begin{minipage}{.46\linewidth}
        \centering
        \includegraphics[scale=0.4]{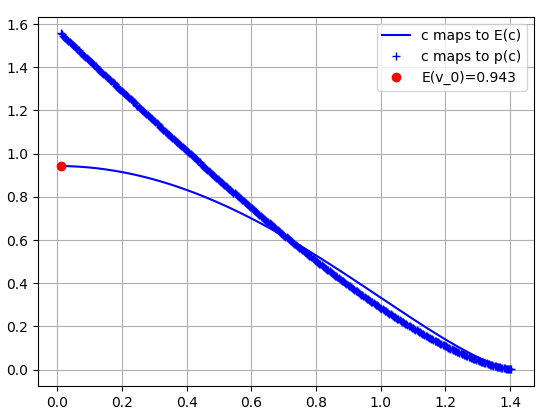}
    \end{minipage}
    \hfill%
    \begin{minipage}{.46\linewidth}
        \centering
        \includegraphics[scale=0.4]{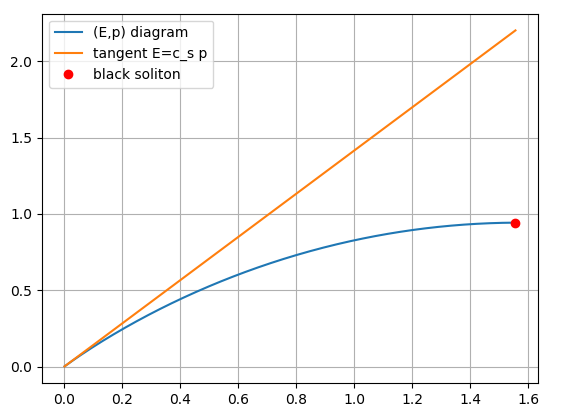}
    \end{minipage}
    \captionof{figure}{Numerical simulations with $f(\rho)=1-\rho$ (Gross-Pitaevskii).}
\end{center}

Now we give two examples of the nonlinearity introduced in Remark~\ref{remarque: exemple de non linearité f(x)=1-x+ a(1-x) ^2p-1} for which there exists a unique travelling wave at any speed $c$. We observe two drastically different behaviours. In the first plot, we have chosen $a$ small compared to $p$, and we have a curve that looks like the Gross-Pitaevskii one.

\begin{center}
    \begin{minipage}[c]{.46\linewidth}
        \centering\includegraphics[scale=0.4]{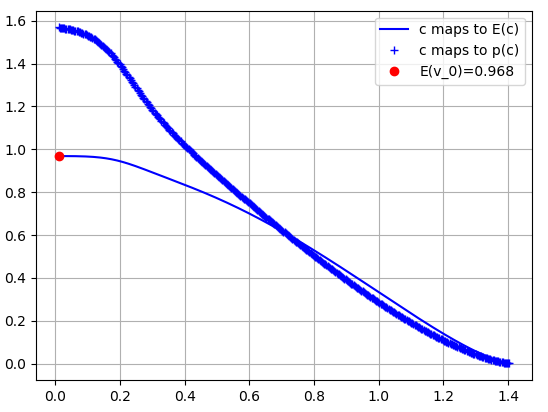}
    \end{minipage}
    \hfill%
    \begin{minipage}[c]{.46\linewidth}
        \centering\includegraphics[scale=0.4]{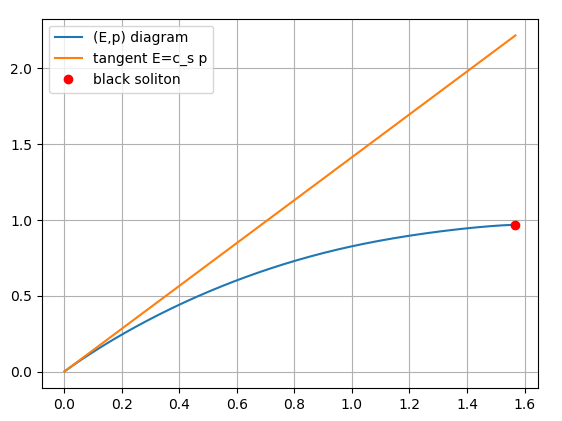}
    \end{minipage}
    \captionof{figure}{Numerical simulations with $f(\rho)=1-\rho + 10(1-\rho)^{59}$.}
\end{center}

On the other hand, taking $a$ widely larger than $p$, we observe a change of variation for the function $c\mapsto p(\gv_c)$, which corresponds to a cusp in the energy/momentum diagram. Therefore, according to Figure~3, there has to be a soliton on the right side of the red dot with the same energy than the black soliton, this is interpreted as a soliton  with a speed $c_*$ such that $p(\gv_{c_*})=\gq_*$ and the red star stands for it.

\begin{center}
    \begin{minipage}[c]{.46\linewidth}
        \centering
        \includegraphics[scale=0.4]{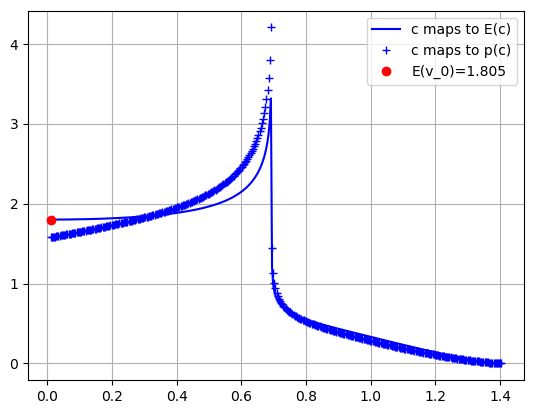}
    \end{minipage}
    \hfill%
    \begin{minipage}[c]{.46\linewidth}
        \centering
        \includegraphics[scale=0.4]{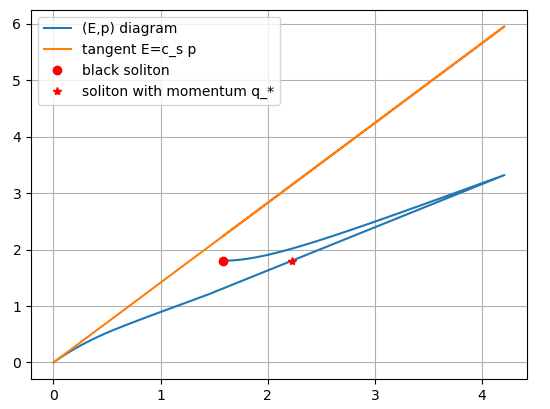}
    \end{minipage}
    \captionof{figure}{Numerical simulations with $f(\rho)=1-\rho + 120(1-\rho)^{19}$.}
\end{center}

Let us now study $f(\rho)=4(1-\rho)+36(1-\rho)^3$. This nonlinearity is widely investigated in Example 2 in~\cite{Chiron7}. The set of speeds where there exists a non trivial travelling wave is $[0,c_0)\cup (c_0,c_s)$ for some $0 < c_0 < c_s$, and there is an asymptote for both branches of the $(E,p)$ diagram when $c\rightarrow c_0^+$ for the lower branch in green (resp. when $c\rightarrow c_0^-$ for the convexe part of the upper branch in blue). This asymptote can be shown to be the line with equation $E=c_0 p + E_0$ where $E_0\approx 0.0512$.

\begin{center}
    \begin{minipage}[c]{.46\linewidth}
        \centering
        \includegraphics[scale=0.4]{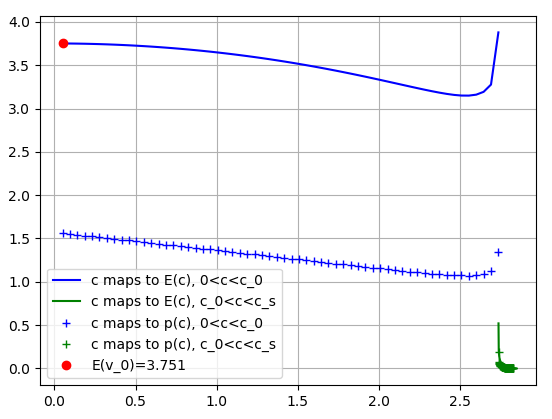}
    \end{minipage}
    \hfill%
    \begin{minipage}[c]{.46\linewidth}
        \centering
        \includegraphics[scale=0.4]{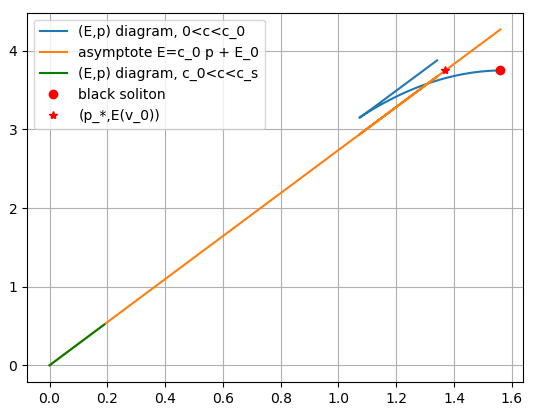}
    \end{minipage}
    \captionof{figure}{Numerical simulations with $f(\rho)=4(1-\rho) + 36(1-\rho)^{3}$.}
\end{center}

There is also a cusp because of the change of variation of the momentum. However, because of the singularity at $c_0$, this gives rise, according to Figure~4, to a concave part in the $(E,p)$ diagram corresponding to small speeds. For numerical reasons, one cannot plot the values of $p(\gv_c)$ when $c$ is to close to $c_0$, this gives a tremendous lack of data regarding the lower branch of the $(E,p)$ diagram. However, since it is supposed to be an asymptote, we can expect the green curve to remain close to the orange tangent as $p$ rises, so that the value of $\gq_*$ is close to the momentum $\gp_*\approx 1.37< \frac{\pi}{2}$ where the tangent passes through the point $\big(\gp_*,E(\gv_0)\big)$ represented by the red star.

\numberwithin{equation}{section}

\section{Appendix}
\appendix

\section{A few estimates related to the momentum and the energy}

A first step to have control on the momentum or on a function in terms of its energy is the following result.

\begin{lem}\label{lemme majoration pointwise de F(|rho|^2)phi' par e(v) su rho}
Let $\rho$ and $\varphi$ be two real-valued, smooth functions on some interval of $\R$. Assume that $\rho$ is positive and set $v=\rho e^{i\varphi}$. Then, we have the pointwise bound

\begin{equation*}\label{lem: majoration pointwise de F(|rho|^2)phi' par e(v) su rho}
    \sqrt{F(\rho^2)}|\varphi'|\leq \dfrac{e(v)}{\rho}.
\end{equation*}
\end{lem}

\begin{proof}
This is a consequence of Young's inequality $2ab\leq a^2+b^2$ on the real line, for $a=\rho\varphi'$ and $b=\sqrt{F(\rho^2)}$, when we decompose the energy as

$$e(v)=\dfrac{1}{2}\Big( \rho'^2+ (\rho\varphi')^2 + \big(\sqrt{F(\rho^2)}\big)^2\Big)\geq \dfrac{1}{2}\big( (\rho\varphi')^2 + \big(\sqrt{F(\rho^2)}\big)^2\big).$$
\end{proof}

Using the previous lemma, we are led to the following pointwise control on a function $v\in\Nenergyset$ with respect to its energy and momentum.

\begin{cor}\label{corollaire du lemme sqrt(f rho^2)|varphi'| leq lambda e(v)/rho}
Consider a function $v=\rho e^{i\varphi}\in \mathcal{N}\energyset$ and assume that \eqref{1ere hypothèse de croissance sur F minorant intermediaire} holds. Then we have 

\begin{equation}\label{cor: 1ere inegalité}
    |(1-\rho^2)\varphi'|\leq \lambda^{-\frac{1}{2}}\dfrac{e(v)}{\rho}.
\end{equation}

In particular, if there exists  $\varepsilon\in (0,1)$ such that $1-\varepsilon\leq |v|^2$ uniformly on $\R$, then

\begin{equation}\label{cor: 2eme inegalité}
    |p(v)|\leq \dfrac{E(v)}{\sqrt{4\lambda(1-\varepsilon)}}.
\end{equation}

For $p(v)=\gp >0$, we get
\begin{equation}\label{cor: 3eme inegalité}
    \inf_\R |v|\leq (4\lambda)^{-\frac{1}{2}}\dfrac{E(v)}{\mathfrak{p}}.
\end{equation}

In particular, if $\delta(v):=1-(4\lambda)^{-\frac{1}{2}}\dfrac{E(v)}{\mathfrak{p}} >0$, then, given any $\delta\in \big(0,\delta(v)\big)$, there exists $x_\delta\in \R$ such that 

\begin{equation}\label{corollaire: 1-|v(x_delta)| geq delta}
    1-|v(x_\delta)| \geq \delta.
\end{equation}

\end{cor}

\begin{proof}
Inequality \eqref{cor: 1ere inegalité} is a straightforward consequence of Lemma \ref{lemme majoration pointwise de F(|rho|^2)phi' par e(v) su rho} combined with the hypothesis~\eqref{1ere hypothèse de croissance sur F minorant intermediaire}. Estimate \eqref{cor: 2eme inegalité} follows integrating \eqref{cor: 1ere inegalité}. Finally, \eqref{cor: 3eme inegalité} is based on the same argument as previously replacing $\sqrt{1-\varepsilon}$ by $\inf_\R |v|$.
\end{proof}

Secondly, we state a result that gives a control of the uniform norm for the function $\eta$ in terms of its energy.

\begin{lem}\label{lemme estimée de la norme infinie de eta en fonction de l'energie E(v)}
Suppose that $F$ verifies \eqref{1ere hypothèse de croissance sur F minorant intermediaire}. Let $v\in \energyset$ and set $\eta:= 1 -|v|^2$. Then

\begin{equation}\label{lem: estimée de la norme infinie de eta en fonction de l'energie E(v)}
    \Vert \eta\Vert^2_{L^\ii}\leq E(v)\big(A+BE(v)\big),
\end{equation}

with $A=\dfrac{64}{\lambda}$ and $B=\dfrac{16}{\sqrt{\lambda}}$.
\end{lem}

\begin{proof}
For $x\in\R$, we define $\rho(x):=|v(x)|$. Since $F$ satisfies~\eqref{1ere hypothèse de croissance sur F minorant intermediaire} and $\eta\underset{\pm \ii}{\rightarrow}0$,

\begin{align*}
    \eta(x)^2 &=\int_{-\ii}^x 2\eta\eta'\\
    &\leq \sqrt{\lambda}\int_\R \eta^2 +\dfrac{1}{\sqrt{\lambda}}\int_\R \eta'^2\\
    &\leq \dfrac{1}{\sqrt{\lambda}}\bigg(\int_\R F(\rho^2) + 4\int_\R \rho^2 (\rho')^2\bigg)\\
    &\leq \dfrac{2}{\sqrt{\lambda}}\Big(E_p(v)+4\Vert v\Vert_{L^\ii}^2 E_k(v)\Big)\\
    &\leq \dfrac{2}{\sqrt{\lambda}}\max\big( 1,4\Vert v\Vert_{L^\ii}^2\big)E(v).
\end{align*}

We deduce that 

\begin{equation*}
    \normLii{\eta}^2\leq \dfrac{8}{\sqrt{\lambda}}(1+\normLii{\eta})E(v).
\end{equation*}

Solving this inequality, we get

\begin{equation}\label{preuve lemme appendix normlii eta inf (a+sqrt a(a+4))/2}
    \normLii{\eta}\leq \dfrac{a + \sqrt{a(a+4)}}{2},
\end{equation}

with $a=\dfrac{8}{\sqrt{\lambda}}E(v)$, so that we obtain

\begin{equation*}
    \normLii{\eta}^2 \leq E(v)\big(A+BE(v)\big)\quad\text{with }A=\dfrac{64}{\lambda}\text{ and }B=\dfrac{16}{\sqrt{\lambda}}.
\end{equation*}

\end{proof}

\section{Convergences}

In this short subsection, we compare the various topologies at hand. For instance, we prove that $\energyset$ is weaker when it is endowed with the $H^1$-convergence\footnote{More precisely, the fact that $v_n\rightarrow v$ for the $H^1$-convergence means that $1-|v_n|^2$ (or $1-|v_n|$) tends to $1-|v|^2$ for the $H^1$-norm.} than the metric space $\big(\energyset,d_A\big)$.

\begin{lem}\label{lemme convergence pour d_A implique la convergence dans H^1 de 1-rho_n^2 vers 1-rho^2}
Suppose that we have a sequence $(v_n)\in\mathcal{N}\energyset^{\N}$ and $v\in\mathcal{N}\energyset$ such that
\begin{equation*}
    v_n\overset{d_{A}}{\ntend} v .
\end{equation*}

Then, we have

\begin{equation}\label{lemme convergence pour d_A implique convergence dans H^1, convergence dans H^1}
    1-|v_n|^2\ntend 1-|v|^2\quad\text{in }H^1(\R).
\end{equation}

In particular, 
\begin{equation*}
    p(v_n)\ntend p(v).
\end{equation*}

\end{lem}

\begin{proof}
The convergence of $1-|v_n|^2$ in $L^2(\R)$ is straightforward. Now, let us prove that 

\begin{equation}\label{preuve de lemme langle v_n,v_n' rangle tend dans L^2 vers}
    \langle v_n,v_n'\rangle_\C\ntend \langle v,v'\rangle_\C\quad\text{in }L^2(\R).
\end{equation}

We have 
\begin{equation*}
    -\langle v_n,v_n'\rangle_\C+\langle v,v'\rangle_\C =\langle v_n,v'-v'_n\rangle_\C + \langle v-v_n,v'\rangle_\C .
\end{equation*}

On the other hand, we have
\begin{equation*}
    \normLdeux{1-|v_n|}\leq\normLdeux{1-|v_n|^2}\quad\text{and}\quad \normLdeux{\partial_x (1-|v_n|)}\leq \normLdeux{\partial_x v_n},
\end{equation*}

so that, by the Sobolev embedding theorem,
\begin{equation*}
    \normLii{1-|v_n|}\leq C\big(\normLdeux{1-|v_n|^2}+\normLdeux{\partial_x v_n}\big).
\end{equation*}

In view of the convergence for $d_A$, there exists some $M>0$ such that for all $n\in\N$,

\begin{equation*}
    \normLii{1-|v_n|}\leq M.
\end{equation*}

We obtain 

\begin{align*}
    \normLdeux{\langle v_n,v'-v'_n\rangle_\C} & \leq \big(\normLii{1-|v_n|}+1)\normLdeux{v'-v_n'}\\
    &\leq (M+1)\normLdeux{v'-v_n'}\ntend 0.
\end{align*}

In another direction, for $\varepsilon >0$, we can exhibit $B_\varepsilon > A >0$ such that 
\begin{equation*}
    \Vert v'\Vert_{L^2([-B_\varepsilon,B_\varepsilon]^c)}\leq\dfrac{\varepsilon}{M+1+\normLii{v}},
\end{equation*} 

because $v'\in L^2(\R)$. We also have the estimate

\begin{align*}
    \normLdeux{\langle v-v_n,v'\rangle_\C} & \leq \Vert v-v_n\Vert_{L^{\ii}([-B_\varepsilon,B_\varepsilon])}\Vert v'\Vert_{L^{2}([-B_\varepsilon,B_\varepsilon])} \\
    & + \big( \Vert v\Vert_{L^{\ii}([-B_\varepsilon,B_\varepsilon]^c)} + \Vert v_n\Vert_{L^{\ii}([-B_\varepsilon,B_\varepsilon]^c)} \big)\Vert v' \Vert_{L^{2}([-B_\varepsilon,B_\varepsilon]^c)}.
\end{align*}

In particular, we get

\begin{align*}
    \big( \Vert v\Vert_{L^{\ii}([-B_\varepsilon,B_\varepsilon]^c)} + \Vert v_n\Vert_{L^{\ii}([-B_\varepsilon,B_\varepsilon]^c)} \big)\Vert v' \Vert_{L^{2}([-B_\varepsilon,B_\varepsilon]^c)} & \leq \big( M+1 + \normLii{v}\big) \Vert v' \Vert_{L^{2}([-B_\varepsilon,B_\varepsilon]^c)} \\
    & \leq\varepsilon .
\end{align*}

To deal with the other term, we write that for all $x\in\R$,

\begin{align*}
    v_n (x)-v(x)=v_n(0)-v(0) +\int_0^x \big( v_n'(t)- v'(t)\big)dt .
\end{align*}

This leads to
\begin{align*}
    \Vert v_n -v\Vert_{L^{\ii}([-B_\varepsilon,B_\varepsilon])} & \leq \Vert v_n -v\Vert_{L^{\ii}([-A,A])} + \sqrt{B_\varepsilon}\normLdeux{v_n'-v'}\\
    & \leq (1+\sqrt{B_\varepsilon}) d_A (v_n,v).
\end{align*}

We take $n$ large enough so that $(1+\sqrt{B_\varepsilon}) d_A (v_n,v)\leq \varepsilon$ and we infer that there exists an integer $N$ such that for all $n\geq N$,

\begin{equation*}
    \normLdeux{\langle v_n,v_n'\rangle_\C-\langle v,v'\rangle_\C}\leq 2\varepsilon .
\end{equation*}

This concludes the proof of \eqref{lemme convergence pour d_A implique convergence dans H^1, convergence dans H^1}.

Furthermore, since $v_n,v\in\Nenergyset$, this allows us to write the liftings $v_n = \rho_n e^{i\varphi_n}$ and $v = \rho e^{i\varphi}$. We have 

\begin{equation*}
    \varphi_n' = \dfrac{\langle iv_n,(v_n)'\rangle}{\rho_n^2}\quad\text{and}\quad \varphi'=\dfrac{\langle iv,v'\rangle}{\rho^2}.
\end{equation*}

In particular, as a consequence of~\eqref{lemme convergence pour d_A implique convergence dans H^1, convergence dans H^1}, the Sobolev embedding theorem and the fact that $v\in\Nenergyset$, one can find $m>0$ independent of $n$ such that 

\begin{equation*}
    \inf_{\R} |v_n| \geq m\quad\text{and}\quad\inf_\R |v|\geq m ,
\end{equation*}

so that we can estimate

\begin{align}
    \normLdeux{\varphi_n ' - \varphi'} & \leq \Big\Vert \dfrac{\rho^2 \langle iv_n,(v_n)'\rangle - \rho_n^2 \langle iv,v'\rangle }{\rho_n^2\rho^2}\Big\Vert_{L^2}\notag\\
    &\leq \dfrac{1}{m^4} \big(\normLii{\rho^2}\normLdeux{\langle iv_n,(v_n)'\rangle - \langle iv,v'\rangle}+\normLdeux{\langle iv,v'\rangle}\normLii{\rho_n^2-\rho^2}\big).\label{varphi_n' -varphi_c' inférieure en norme L2}
\end{align}

In view of \eqref{lemme convergence pour d_A implique convergence dans H^1, convergence dans H^1}, we deduce that the right-hand term in \eqref{varphi_n' -varphi_c' inférieure en norme L2} tends to $0$ (we refer to the same calculus as previously for the convergence of $\normLdeux{\langle iv_n,(v_n)'\rangle - \langle iv,v'\rangle}$ to 0). We infer that

\begin{equation*}
    \varphi_n'\ntend\varphi'\quad\text{and}\quad 1-\rho_n^2\ntend 1-\rho^2\quad\text{in }L^2(\R),
\end{equation*}

and therefore
\begin{equation*}
    p(v_n)=\dfrac{1}{2}\int_\R (1-\rho_n^2)\varphi_n' \ntend \dfrac{1}{2}\int_\R (1-\rho^2)\varphi'=p(v).
\end{equation*}

This finishes the proof of Lemma~\ref{lemme convergence pour d_A implique la convergence dans H^1 de 1-rho_n^2 vers 1-rho^2}.

\end{proof}

The untwisted momentum is also continuous with respect to the distance $d_A$, we refer to~\cite{BeGrSaS1} for the proof of this lemma.

\begin{lem}[\cite{BeGrSaS1}]\label{lemme appendix untwisted momentum continu pour d_A}
    Suppose that we have a sequence $(v_n)\in\energyset^{\N}$ and $v\in\energyset$ such that
\begin{equation*}
    v_n\overset{d_{A}}{\ntend} v .
\end{equation*}

Then,
\begin{equation*}
    [p](v_n)\ntend [p](v).
\end{equation*}
\end{lem}

We now deal with the convergence of the potential term of the energy.

\begin{lem}\label{lemme convergence de 1-rho_n dans H^1 implique convergence de F(rho_n^2)}
Suppose that \eqref{hypothèse de croissance sur F majorant} holds. Consider a sequence $(\rho_n)\in 1+ H^1(\R)^{\N}$ and $\rho \in 1+ H^1(\R)$ such that
\begin{equation*}
    1-\rho_n\ntend 1-\rho\quad\text{in }H^1(\R).
\end{equation*}

Then

\begin{equation}
    F(\rho_n^2)\ntend F(\rho^2)\quad\text{in }L^1(\R).
\end{equation}
\end{lem}

\begin{proof}
By the Sobolev embedding theorem, we have 

\begin{equation}\label{convergence de vert rho_n - rho vert_p pour tout p entre 2 et +ii}
    \Vert \rho_n -\rho \Vert_{L^p}\ntend 0 \quad\text{for all }p\in [2,+\ii].
\end{equation}

Since $1-\rho_n^2 - (1- \rho^2 ) =(\rho_n - \rho) ( 1- \rho_n +1- \rho - 2)$, we get 
\begin{equation*}
    \Vert 1-\rho_n^2 -(1- \rho^2) \Vert_{L^\ii}\leq \Vert \rho_n - \rho \Vert_{L^p} (\Vert 1-\rho_n \Vert_{L^\ii} +\Vert 1-\rho\Vert_{L^{\ii}} +2 )\ntend 0\quad\text{for all }p\in [2,+\ii].
\end{equation*}

Provided that $F(1)=F'(1)=0$, the Taylor formula on $[0,2]$ yields a constant $C\geq 0$ such that

\begin{equation*}
    F(\xi)\leq C(\xi - 1 )^2\quad\text{for  }\xi\in [0,2].
\end{equation*}

The previous inequality, combined with assumption \eqref{hypothèse de croissance sur F majorant} implies that for all $\xi\in\R_+$,

\begin{equation*}
    F(\xi)\leq C(\xi - 1)^2 + M|\xi - 1 |^q ,
\end{equation*}

so that 

\begin{equation*}
    F(\rho^2)\leq C(1-\rho^2)^2 + M|1-\rho^2|^q .
\end{equation*}

We show that, up to a subsequence, there exist $g_1\in L^2(\R)$ and $g_2 \in L^q(\R)$ such that 

\begin{equation*}
    \left\{
\begin{array}{l}
    1- \rho_n^2 \ntend 1-\rho^2\quad a.e. \\
    |1-\rho_n^2|\leq g_1 \quad a.e. \\
    |1-\rho_n^2|\leq g_2 \quad a.e. 
    \end{array}
\right.
\end{equation*}

Indeed, we already showed that \begin{equation*}
    1-\rho_n^2 \ntend 1-\rho^2\quad\text{in }L^2(\R),
\end{equation*}

then, up to a subsequence, there exists $g_1 \in L^2$, such that 

\begin{equation*}
    \left\{
\begin{array}{l}
    1- \rho_n^2 \ntend 1-\rho^2\quad a.e. \\
    |1-\rho_n^2|\leq g_1 \quad a.e.
    \end{array}
\right.
\end{equation*}

It remains to verify that the same property holds in $L^q(\R)$. This is true because of \eqref{convergence de vert rho_n - rho vert_p pour tout p entre 2 et +ii} and the fact that $1-\rho_n^2 \in H^1(\R)$ for all $n\in\N$.

Therefore, we deduce from the continuity of $F$ that 

\begin{equation*}
    \left\{
\begin{array}{l}
    F(\rho_n^2) \ntend F(\rho^2) \quad a.e. \\
    F(\rho_n^2) \leq  C g_1^2 + M g_2^q \quad a.e.  \\
    \end{array}
\right.
\end{equation*}

By Lebesgue's dominated convergence theorem, we are led to

\begin{equation*}
    F(\rho_n^2)\ntend F(\rho^2)\quad\text{in }L^1(\R).
\end{equation*}

For any subsequence of $\big(F(\rho_n^2)\big)$, the same argument holds and we can always find a further subsequence that converges to $F(\rho^2)$. Therefore we get the complete convergence in $L^1(\R)$.

\end{proof}

\section{Localisation of $\gq_*$}\label{subsection localisation de q_*}

Whenever there exists a complete branch $c\mapsto \gv_c\in \mathcal{C}^1\big((0,c_s),\Nenergyset\big)$, then $c\mapsto p(\gv_c)$ is $\mathcal{C}^1$on $(0,c_s)$ (because $\gv_c$ does not vanish on $(0,c_s)$). Moreover, by Proposition~\ref{prop: la solution de tw0 est unique si int_0^1 f neq 0}, the kink exists a priori according to D. Chiron's work (Lemma~4 in~\cite{Chiron8}) we have the limits $p(\gv_c)\rightarrow \frac{\pi}{2}$ and $E(\gv_c)\rightarrow E(\gv_0)$ as $c$ tends to $0$. In this framework, we show that if the assumptions of Theorem~\ref{existence et unicité des travelling wave si condition du zéro} are achieved, and if moreover there exists a finite number of speeds such that $\frac{d}{dc}p(\gv_c)$ vanishes, then $\frac{\pi}{2}\leq \gq_*$.

We split the argument in two parts
and we first assume that $\frac{d}{dc}\big(p(\gv_c)\big)_{|c=c_*}=0$ for a unique $c_*\in (0,c_s)$. Set $\gp_*:=p(\gv_{c_*})$ and consider the set \begin{equation}\label{definition courbe énergie moment}
    \mathcal{G}:=\big\lbrace \big(p(\gv_c),E(\gv_c)\big)\big|c\in (0,c_s)\big\rbrace .
\end{equation}

The function $c\mapsto p(\gv_c)$ is necessarily decreasing on $(c_*,c_s)$. Indeed, according to the inverse function theorem, there exist two branches $\mathcal{E}^{CCV}$ and $\mathcal{E}^{CVX}$ whose union is equal to $\mathcal{G}$. By uniqueness and Theorem~\ref{theoreme existence de travelling wave}, $\Emin = \mathcal{E}^{CCV}$ on $(0,\gq_*)$ and $\Emin$ is concave on $(0,\gq_*)$. In view of~\eqref{le signe de mathcal E'' est le meme que celui de p'}, there is no other possibility that $\mathcal{E}^{CCV}$ is strictly concave on $(0,\gp_*)$ and that $\mathcal{E}^{CVX}$ is convex on $(\frac{\pi}{2},\gp_*)$. In the case where $c\mapsto p(\gv_c)$ is also decreasing on $(0,\gp_*)$, then $c\mapsto p(\gv_c)$ performs a global diffeomorphism between $(0,c_s)$ and $(0,\frac{\pi}{2})$ and it implies that $\gq_*=\frac{\pi}{2}$. 

Suppose by contradiction that $\gq_* <\frac{\pi}{2}$. In view of the previous work, we have the following variations.
$$\begin{tikzpicture}
   \tkzTab{$c$ / 1 , $\dfrac{d}{dc}p(\gv_c)$ / 1, $p(\gv_c)$ / 1.5}{$0$, $c_*$, $c_s$}
          {, +, z, -, }
          {-/ $\frac{\pi}{2}$, +/ $\gp_*$, -/ $0$}
\end{tikzpicture}
$$

By concavity (resp. convexity), $\mathcal{E}^{CCV}$ (resp. $\mathcal{E}^{CVX}$) always lies under (resp. over) its tangent at the point $\gp_*$. Using this global inequality at the point $\gq_*$ (resp. $\frac{\pi}{2})$, we derive
\begin{align*}
    E(\gv_0)=\mathcal{E}^{CCV}(\gq_*)&\leq (\mathcal{E}^{CCV})'(\gp_*)(\gq_*-\gp_*)+\mathcal{E}^{CCV}(\gp_*)=c_* (\gq_*-\gp_*)+E(\gv_{c_*})\\
    & < c_* \Big(\frac{\pi}{2}-\gp_*\Big) + E(\gv_{c_*})\leq \mathcal{E}^{CVX}\Big(\frac{\pi}{2}\Big)=E(\gv_0),
\end{align*}

which brings a contradiction.

Now, we generalize this argument when $c\mapsto \frac{d}{dc}\big(p(\gv_c)\big)$ vanishes at the speeds $0<c_1 <...< c_J<c_s$. We assume by contradiction that $\gq_* <\frac{\pi}{2}$, then we first have $\frac{d}{dc}\big(p(\gv_c)\big)>0$ on $(0,c_1)$. Otherwise, we would also have $\frac{d}{dc}\big(E(\gv_c)\big)<0$ by~\eqref{relation de groupe hamiltonien dE/dc(v_c)=cdp/dc(v_c)}. Now taking $c\in (0,c_1)$ small enough so that $\gq_*<p(\gv_c)<\frac{\pi}{2}$ and using Proposition~\ref{prop Emin(q) leq Emin(pi/2) = E(v_0)}, we obtain the following contradiction
\begin{equation*}
    E(\gv_0)=\Emin(\gq_*)=\Emin\big(p(\gv_c)\big)\leq E(\gv_c)<E(\gv_0).
\end{equation*}

This allows us to define
\begin{equation*}
    j_*:=\min\Big\{ j\in\{ 1,...,J\}\Big| \dfrac{d}{dc}\big(p(\gv_c)\big)>0\text{ on }(c_{j-1},c_j)\text{ and } \frac{d}{dc}\big(p(\gv_c)\big)<0\text{ on }(c_j,c_{j+1})\Big\},
\end{equation*}

with the convention $(c_0,c_{J+1}):=(0,c_s)$. We argue as above, we write $\mathcal{E}^{CVX}$ (resp. $\mathcal{E}^{CCV}$) the branch of the energy/momentum defined by the graph~\eqref{definition courbe énergie moment} for speeds in $(0,c_{j_*})\setminus \{c_1,...,c_{j_*-1}\}$ (resp. for speeds in $(c_{j_*},c_{j_*+1})$). Furthermore, by definition of $j_*$, we can extend both these functions as follows: $\mathcal{E}^{CVX}\big(p(\gv_{c_{k}})\big)=E(\gv_{c_{k}})$ for $k\in\{1,...,j_*-1\}$, $\mathcal{E}^{CCV}\big(p(\gv_{c_{j_*+1}})\big)=E(\gv_{c_{j_*+1}})$ and $\mathcal{E}^{CCV}\big(p(\gv_{c_{j_*}})\big)=\mathcal{E}^{CVX}\big(p(\gv_{c_{j_*}})\big)=E(\gv_{c_{j_*}})$. By concavity of $\mathcal{E}^{CCV}$ (resp. convexity of $\mathcal{E}^{CVX}$), we deduce as above that
\begin{align*}\label{mathcal{E}^{CCV}big(p(gv_{c_{j_*+1}})big)eq c_{j_*}}
    \mathcal{E}^{CCV}\big(p(\gv_{c_{j_*+1}})\big)\leq c_{j_*} \big( p(\gv_{c_{j_*+1}})- p(\gv_{c_{j_*}})\big) + \mathcal{E}\big(p(\gv_{c_{j_*}})\big)< \mathcal{E}^{CVX}\Big(\dfrac{\pi}{2}\Big)=E(\gv_0).
\end{align*}

If $E(\gv_0)\leq \mathcal{E}^{CCV}\big(p(\gv_{c_{j_*+1}})\big)$, the inequality above brings a contradiction. Otherwise, $E(\gv_0)> \mathcal{E}^{CCV}\big(p(\gv_{c_{j_*+1}})\big)=E(\gv_{c_{j_*+1}})$, then by Lemma~\ref{lemme l'inf de l'énergie dont la fonction s'annule est atteint en u_0}, we must have $p(\gv_{c_{j_*+1}})\leq \gq_* <\frac{\pi}{2}$, then $\mathcal{E}^{CCV}$ is well-defined on $\gq_*$ and we have \begin{equation*}\label{mathcal{E}^{CCV}big(p(gv_{c_{j_*+1}})big)eq c_{j_*}}
    E(\gv_0)=\Emin(\gq_*)\leq \mathcal{E}^{CCV}(\gq_*)\leq c_{j_*} \big( \gq_*- p(\gv_{c_{j_*}})\big) + \mathcal{E}\big(p(\gv_{c_{j_*}})\big)< \mathcal{E}^{CVX}\Big(\dfrac{\pi}{2}\Big)=E(\gv_0),
\end{equation*} 
which provides an ultimate contradiction.

\subsection*{Acknowledgements.}
I am grateful to the reviewers for considering my article and for their advice and contributions. I am also thankful to P. Gravejat for his caring support over the past two years. This work was supported by the CY Initiative of Excellence (Grant “Investissements d’Avenir” ANR-16-IDEX-0008).

\bibliographystyle{plain}
\bibliography{Bibliogr}

\end{document}